\documentclass[11pt,leqno]{article}
\usepackage{a4wide}
\usepackage{amssymb}
\usepackage{amsmath}
\usepackage{amsthm}
\usepackage{mathrsfs}
\usepackage{yfonts}
\usepackage{euscript}
\usepackage{upgreek}
\usepackage[all]{xypic}
\usepackage[usenames]{color}
\usepackage[usenames,dvipsnames]{xcolor} %{\color{Maroon/ForestGreen}  }
\usepackage{mathabx} % f\"{u}r \corresponds = entspricht-Zeichen, und f\"ur \varhash
\usepackage[ngerman,shadow,color=blue!30,textsize=footnotesize,colorinlistoftodos]{todonotes}
\usepackage[T1]{fontenc}
\usepackage[utf8]{inputenc}
\usepackage{comment}
\usepackage{xr-hyper}
\usepackage[unicode,psdextra]{hyperref}
\usepackage{bookmark}
\usepackage{enotez, translations}

\newcommand{\Footnote}[1]{\footnote{{\color{Sepia} #1}}}  % Fussnoten nur f\"{u}r Peter und Otmar farbig in Sepia
\renewcommand{\Footnote}[2][]{\relax}  % nur Fussnoten f\"{u}r Peter und Otmar ausblenden
\theoremstyle{plain}

\newcommand{\id}{\operatorname{id}}

\newcommand{\Hom}{\operatorname{Hom}}

\newcommand{\Ext}{\operatorname{Ext}}

\newcommand{\Rep}{\operatorname{Rep}}

\newcommand{\Aut}{\operatorname{Aut}}

\newcommand{\Gal}{\operatorname{Gal}}
\newcommand{\rk}{\operatorname{rk}}

 \newcommand{\QQ}{\mathbb{Q}}

 \newcommand{\cC}{\mathcal{C}}

\newcommand{\cT}{\mathcal{T}}

\renewcommand{\cR}{\mathcal{R}}

 \newcommand{\Qp}{\QQ_p}

 \newcommand{\be}{\begin{equation}}
\newcommand{\ee}{\end{equation}}

\newcommand{\fpbar}{\ifmmode {\overline{\mathbb{F}_p}}\else$\mathbb{F}_p$\ \fi}

\newcommand{\fp}{\ifmmode {\mathbb{F}_p}\else$\mathbb{F}_p$\ \fi}
\newcommand{\zp}{\ifmmode \mathbb{Z}_p\else$\mathbb{Z}_p$\ \fi}
\newcommand{\zpur}{\ifmmode \widehat{\zp^{ur}}\else $\widehat{\zp^{ur}}$\ \fi}

 \newcommand{\TLT}{T_{\pi}}

\newcommand{\La}{\ifmmode\Lambda\else$\Lambda$\fi}

\newcommand{\q}{\ifmmode {\mathbb Q}\else${\mathbb Q}$\ \fi}
\newcommand{\qp}{\ifmmode {\mathbb Q}_p\else${\mathbb Q}_p$\ \fi}
\newcommand{\z}{\mathbb{Z}}

\newcommand{\Q}{\ifmmode {\mathbb Q}\else${\mathbb Q}$\ \fi}
\newcommand{\ql}{\ifmmode {{\mathbb Q}_l}\else${\mathbb Q}_l$\ \fi}

\newtheorem{theorem}{Theorem}[section]
\newtheorem{corollary}[theorem]{Corollary}
\newtheorem{lemma}[theorem]{Lemma}
\newtheorem{remark}[theorem]{Remark}
\newtheorem{proposition}[theorem]{Proposition}

\newtheorem{definition}[theorem]{Definition}

\theoremstyle{remark}

\author{Peter Schneider and Otmar Venjakob}
\begin{document}

\title{\textbf{Compairing categories of Lubin-Tate $(\varphi_L,\Gamma_L)$-modules}}

\maketitle

\begin{abstract}
 In the Lubin-Tate setting we compare different categories of $(\varphi_L,\Gamma)$-modules over various perfect or imperfect coefficient rings. Moreover, we  study their associated  Herr-complexes. Finally, we show that a Lubin Tate extension gives rise to a weakly decompleting, but not decompleting tower in the sense of Kedlaya and Liu.
\end{abstract}

\tableofcontents

\section{Introduction}

Since its invention by Fontaine in \cite{Fo} the concept of $(\varphi,\Gamma)$-modules (for the $p$-cyclotomic extension) has become a powerful tool in the study of $p$-adic Galois representations of local fields. In particular, it could be fruitfully applied in Iwasawa theory \cite{benois,B,Na,NaANT,NaRk,V-Kato,LVZ15,LLZ11,BV} and in  the $p$-adic local Langlands programme \cite{Co}. A good introduction to the subject regarding the state of the art around 2010 can be found in \cite{BC, FO}.

Afterwards a couple of generalisations have been developed. Firstly, Berger and Colmez \cite{BeCo}  as well as Kedlaya, Pottharst and Xiao  \cite{KPX} extended the theory to (arithmetic) families of $(\varphi,\Gamma)$-modules, in which representations of the absolute Galois group of a local field on modules over affinoid algebras over $\Qp$ instead of finite dimensional vector spaces are studied. Secondly, parallel to and influenced by Scholze's point of view of perfectoid spaces as well as the upcoming of the Fargues-Fontaine curve \cite{FF} Kedlaya and Liu developed  a (geometric) relative $p$-adic Hodge theory \cite{KLI, KLII}, in which  the  Galois group of a local field is replaced by the \'{e}tale fundamental group of affinoid spaces over $\Qp$ thereby extending an earlier approach by Andreatta and Brinon. In particular, Kedlaya and Liu have introduced systematically $(\varphi,\Gamma)$-moduels over {\it perfect} coefficient rings, i.e., for which the Frobenius endomorphism is surjective, and they have studied  their decent to {\it imperfect} coefficient rings, which is needed for Iwasawa theoretic applications and which generalized the work of Cherbonnier and Colmez \cite{ChCo}.

Recently there has been a growing interest and activity in introducing and studying $(\varphi_L,\Gamma_L)$-modules for Lubin-Tate extensions of a finite extension $L$ of $\Qp$, motivated again by requirements from or potential applications to the  $p$-adic
local Langlands programme \cite{FX,BSX,Co2} or Iwasawa theory \cite{SV15,BF,SV23,MSVW,Poy}. The textbook \cite{GAL}  contains a very detailed and thorough approach to the analogue of Fontaine's original equivalence of categories between Galois representations and \'{e}tale $(\varphi,\Gamma)$-modules to the case of Lubin-Tate extensions as had been proposed, but only sketched in \cite{KR}, see Theorem \ref{KR-equiv}. In this setting it has been shown in \cite{Ku,KV} that - as in the cyclotomic case due to Herr \cite{Her1} - the Galois cohomology of a $L$-representation $V$ of  the absolute Galois group $G_L$ of $L$ can again be obtained as cohomology of a generalized Herr complex for the   $(\varphi_L,\Gamma_L)$-module attached to $V$, see Theorem \ref{thm:Herr}.

The purpose of this article is to spell out in the Lubin-Tate case concretely  the various categories of (classical) $(\varphi_L,\Gamma_L)$-modules over perfect and imperfect coefficient rings (analogously to those considered in \cite{KLI,KLII} who do not cover the Lubin-Tate situation) such as   ${\mathbf{A}}_L, {\mathbf{A}}_L^\dagger, \tilde{\mathbf{A}}_L, \tilde{\mathbf{A}}_L^\dagger, \mathbf{B}_L, {\mathbf{B}}_L^\dagger, {\mathbf{B}}_L, \tilde{\mathbf{B}}_L^\dagger, {\mathcal{R}}_L, \tilde{\mathcal{R}}_L $   to be defined in the course of the main text and to compare them among each other. Moreover, we investigate for which versions the generalized Herr complex calculates again the Galois cohomology of a given representation. The results are summarized in diagrams \eqref{diag1} and \eqref{diag2}. Finally, we study in the last section how Lubin-Tate extensions fit into Kedlaya's and Liu's concept of {\it (weakly) decompleting}
 towers. We show that for $L\neq\Qp$ they are weakly decompleting, but not decompleting.

See \cite{Stein} for some results regarding arithmetic families of  $(\varphi_L,\Gamma_L)$-modules in the Lubin-Tate setting.

{\bf Acknowledgements:}  Both authors are grateful to UBC and PIMS at Vancouver for supporting a fruitful stay.
The project was funded by the Deutsche Forschungsgemeinschaft (DFG, German Research Foundation) – Project-ID 427320536 – SFB 1442, as well as under Germany’s Excellence Strategy EXC 2044 390685587, {\it  Mathematics M\"{u}nster: Dynamics–Geometry–Structure.} We also acknowledge funding by the Deutsche Forschungsgemeinschaft (DFG, German Research Foundation)  under TRR 326 {\it Geometry and
Arithmetic of Uniformized Structures}, project number 444845124, as well as under DFG-Forscher\-grup\-pe award number [1920] {\it Symmetrie, Geometrie und Arithmetik}.

\section{Notation}\label{sec:notation}

Let $\mathbb{Q}_p \subseteq L \subset \mathbb{C}_p$ be a field of finite degree $d$ over $\mathbb{Q}_p$, $o_L$ the ring of integers of $L$, $\pi_L \in o_L$ a fixed prime element, $k_L = o_L/\pi_L o_L$ the residue field,  $q := |k_L|$ and $e$ the absolute ramification index of $L$. We always use the absolute value $|\ |$ on $\mathbb{C}_p$ which is normalized by $|\pi_L| = q^{-1}$. We \textbf{warn} the reader, though, that we will  use the references    \cite{FX} and \cite{La} in which the absolute value is normalized differently from this paper by $|p| = p^{-1}$. Our absolute value is the $d$th power of the one in these references. The transcription of certain formulas to our convention will usually be done silently.

We fix a Lubin-Tate formal $o_L$-module $LT = LT_{\pi_L}$ over $o_L$ corresponding to the prime element $\pi_L$. We always identify $LT$ with the open unit disk around zero, which gives us a global coordinate $Z$ on $LT$. The $o_L$-action then is given by formal power series $[a](Z) \in o_L[[Z]]$. For simplicity the formal group law will be denoted by $+_{LT}$.

%The power series $\frac{\partial (X +_{LT} Y)}{\partial Y}_{|(X,Y) = (Z,0)}$ is a unit in $o_L[[Z]]$ and we let $g_{LT}(Z)$ denote its inverse. Then $g_{LT}(Z) dZ$ is, up to scalars, the unique invariant differential form on $LT$ (\cite{Haz} \S5.8). We also let \begin{equation}\label{f:tLT}
%  \log_{LT}(Z) = Z + \ldots
%\end{equation}
%denote the unique formal power series in $L[[Z]]$ whose formal derivative is $g_{LT}$. This $\log_{LT}$ is the logarithm of $LT$ (\cite{Lan} 8.6). In particular, $g_{LT}dZ = d\log_{LT}$. The invariant derivation $\partial_\mathrm{inv}$ corresponding to the form $d\log_{LT}$ is determined by
%\begin{equation*}
%  f' dZ = df = \partial_\mathrm{inv} (f) d\log_{LT} = \partial_\mathrm{inv} (f) g_{LT} dZ
%\end{equation*}
%and hence is given by
%\begin{equation}\label{f:inv}
%  \partial_\mathrm{inv}(f) = g_{LT}^{-1} f' \ .
%\end{equation}
%For any $a \in o_L$ we have
%\begin{equation}\label{f:dlog}
%  \log_{LT} ([a](Z)) = a \cdot \log_{LT} \qquad\text{and hence}\qquad ag_{LT}(Z) = g_{LT}([a](Z))\cdot [a]'(Z)
%\end{equation}
%(\cite{Lan} 8.6 Lemma 2).

Let $\TLT$ be the Tate module of $LT$. Then $\TLT$ is
a free $o_L$-module of rank one, say with generator $\eta$, and the action of
$G_L := \Gal(\overline{L}/L)$ on $\TLT$ is given by a continuous character $\chi_{LT} :
 G_L \longrightarrow o_L^\times$.

For $n \geq 0$ we let $L_n/L$ denote the extension (in $\mathbb{C}_p$) generated by the $\pi_L^n$-torsion points of $LT$, and we put $L_\infty := \bigcup_n L_n$. The extension $L_\infty/L$ is Galois. We let $\Gamma_L := \Gal(L_\infty/L)$ and $H_L := \Gal(\overline{L}/L_\infty)$. The Lubin-Tate character $\chi_{LT}$ induces an isomorphism $\Gamma_L \xrightarrow{\cong} o_L^\times$.

Henceforth   we use the same notation as in \cite{SV15}.
In particular, the ring endomorphisms induced by sending $Z$ to $[\pi_L](Z)$ are called $\varphi_L$ where applicable; e.g.\  for the ring $\mathscr{A}_L$ defined to be the $\pi_L$-adic completion of $o_L[[Z]][Z^{-1}]$ or $\mathscr{B}_L := \mathscr{A}_L[\pi_L^{-1}]$ which denotes the field of fractions of $\mathscr{A}_L$.  Recall that we also have introduced the unique additive endomorphism $\psi_L$ of $\mathscr{B}_L$ (and then  $\mathscr{A}_L$) which satisfies
\begin{equation*}
  \varphi_L \circ \psi_L = \pi_L^{-1} \cdot trace_{\mathscr{B}_L/\varphi_L(\mathscr{B}_L)} \ .
\end{equation*}
Moreover,  projection formula
\begin{equation*}
  \psi_L(\varphi_L(f_1)f_2) = f_1 \psi_L(f_2) \qquad\text{for any $f_i \in \mathscr{B}_L$}
\end{equation*}
as well as the formula
\begin{equation*}
  \psi_L \circ \varphi_L = \frac{q}{\pi_L} \cdot \id \
\end{equation*}
hold. An  \'{e}tale $(\varphi_L,\Gamma_L)$-module $M$ comes with a Frobenius operator $\varphi_M$ and an induced operator  denoted by $\psi_M$.

Let $\widetilde{\mathbf{E}}^+ := \varprojlim o_{\mathbb{C}_p}/p o_{\mathbb{C}_p}$ with the transition maps being given by the Frobenius $\varphi(a) = a^p$. We may also identify $\widetilde{\mathbf{E}}^+$ with
$\varprojlim o_{\mathbb{C}_p}/\pi_L o_{\mathbb{C}_p}$ with
the transition maps being given by the $q$-Frobenius
$\varphi_q (a) = a^q$. Recall that $\widetilde{\mathbf{E}}^+$ is a complete valuation ring with residue field $\overline{\mathbb{F}_p}$ and its field of fractions $\widetilde{\mathbf{E}} = \varprojlim \mathbb{C}_p$ being algebraically closed of characteristic $p$. Let $\mathfrak{m}_{\widetilde{\mathbf{E}}}$ denote the maximal ideal in $\widetilde{\mathbf{E}}^+$.

The $q$-Frobenius $\varphi_q$ first extends by functoriality to the rings of the Witt vectors $ W(\widetilde{\mathbf{E}})$ and then $o_L$-linearly to $ W(\widetilde{\mathbf{E}})_L := W(\widetilde{\mathbf{E}}) \otimes_{o_{L_0}} o_L$, where $L_0$ is the maximal unramified subextension of $L$. The Galois group $G_L$ obviously acts on $\widetilde{\mathbf{E}}$ and $W(\widetilde{\mathbf{E}})_L$ by automorphisms commuting with $\varphi_q$.  This $G_L$-action is continuous for the weak topology on $W(\widetilde{\mathbf{E}})_L$ (cf.\ \cite[ Lemma 1.5.3]{GAL}).

By sending the variable $Z$ to $\omega_{LT}\in W(\widetilde{\mathbf{E}})_L$ (see directly after \cite[Lem.\ 4.1]{SV15}) we obtain an $G_L$-equivariant, Frobenius compatible  embedding of rings
\begin{equation*}
     \mathscr{A}_L \longrightarrow W(\widetilde{\mathbf{E}})_L \
\end{equation*}
the image of which we call $\mathbf{A}_L$.   The latter ring is a complete discrete valuation ring with prime element $\pi_L$ and residue field the image $\mathbf{E}_L$ of $k_L ((Z)) \hookrightarrow \widetilde{\mathbf{E}}$ sending $Z$ to $\omega:=\omega_{LT} \mod \pi_L.$
We form the maximal integral unramified extension ($=$ strict Henselization) $\mathbf{A}_L^{nr}$  of $\mathbf{A}_L$  inside $W(\widetilde{\mathbf{E}})_L$. Its $p$-adic completion $\mathbf{A}$ still is contained in $W(\widetilde{\mathbf{E}})_L$. Note that $\mathbf{A}$ is a complete discrete valuation ring with prime element $\pi_L$ and residue field the separable algebraic closure $\mathbf{E}_L^{sep}$ of $\mathbf{E}_L$ in $\widetilde{\mathbf{E}}$.  By the functoriality properties of strict Henselizations the $q$-Frobenius $\varphi_q$ preserves $\mathbf{A}$. According to \cite[ Lemma 1.4]{KR} the $G_L$-action on $W(\widetilde{\mathbf{E}})_L$ respects $\mathbf{A}$ and induces an isomorphism $H_L = \ker(\chi_{LT}) \xrightarrow{\cong} \Aut^{cont}(\mathbf{A}/\mathbf{A}_L)$.

Sometimes we omit the index $q, L,$ or $M$ from the Frobenius operator.

Finally, for a valued field K we denote as usual by $\hat{K}$ its completion.
{\color{red}

}

%\section{Perfect and imperfect \texorpdfstring{$(\varphi_L,\Gamma_L)$}{(phi,Gamma)}-modules and their cohomology}

\section{An analogue of Tate's result}

Let $\mathbb{C}_p^\flat $ together with its absolute value $|\cdot|_\flat$ be the tilt of $\mathbb{C}_p $. The aim of this section is to prove an analogue of Tate's classical result \cite[Prop.\ 10]{Ta} for  $\mathbb{C}_p^\flat $ instead of $\mathbb{C}_p $ itself and in the Lubin Tate situation instead of the cyclotomic one. In the following we always consider {\it continuous} group cohomology.

\begin{proposition}\label{prop:TatevanishingH1}
$H^n(H,\mathbb{C}_p^\flat)=0$ for all $n\geq 1$ and $H\subseteq H_L$  any closed subgroup.
\end{proposition}

Since the proof is formally very similar to that of loc.\ cit.\ or  \cite[Prop.\ 14.3.2.]{BC} we only sketch the main ingredients. To this aim we fix $H$ and write sometimes $W$ for $ \mathbb{C}_p^\flat$ as well as $W_{\geq m}:=\{x\in W| |x|_\flat \leq \frac{1}{p^m}\}.$

\begin{lemma}\label{lem:TS1}
The Tate-Sen axiom {\bf (TS1)} is satisfied for $\mathbb{C}_p^\flat$ with regard to $H$, i.e., there exists a real constant $c>1$ such that for all open subgroups $H_1\subseteq H_2$ in $H$ there exists $\alpha\in (\mathbb{C}_p^\flat)^{H_1}$ with $|\alpha|_\flat < c$ and $Tr_{H_2|H_1}(\alpha):=\sum_{\tau\in H_2|H_1 }\tau(\alpha)=1.$ Moreover, for any sequence $(H_m)_m$ of open subgroups $H_{m+1}\subseteq H_m$ of $H$ there exists a trace compatible system $(y_{H_m})_m$ of elements $y_{H_m}\in (\mathbb{C}_p^\flat)^{H_m}$ with $|y_{H_m}|_\flat < c$  and $Tr_{H|H_m}(y_{H_m} )=1$.
\end{lemma}

\begin{proof}
  Note that for a perfect   field $K$ (like $(\mathbb{C}_p^\flat)^{H}$) of characteristic $p$  complete for a multiplicative norm with maximal ideal $\mathfrak{m}_K$ and a  finite extension $F$ one has $\mathrm{Tr}_{F/K}(\mathfrak{m}_F)=\mathfrak{m}_K$ by \cite[Thm.\ 1.6.4]{KedNew}. Fix some $x\in (\mathbb{C}_p^\flat)^{H}$ with $0<|x|_\flat<1$ and set $c:=|x|_\flat^{-1}>1.$ Then we find $\tilde{\alpha}$ in the maximal ideal of $  (\mathbb{C}_p^\flat)^{H_1}$ with $Tr_{H|H_1}(\tilde{\alpha})=x$ and $\alpha:=(Tr_{H_2|H_1}(\tilde{\alpha}))^{-1}\tilde{\alpha}$ satisfies the requirement as $|Tr_{H_2|H_1}(\tilde{\alpha})|_\flat^{-1}\leq |x|_\flat^{-1}=c.$

  For the second claim we successively  choose elements $\tilde{\alpha}_m$ in the maximal ideal of $  (\mathbb{C}_p^\flat)^{H_m}$ such that $Tr_{H|H_1}(\tilde{\alpha}_1)=x$ and $Tr_{H_{m+1}|H_m}(\tilde{\alpha}_{m+1})=\tilde{\alpha}_{m}$ for all $m\geq 1.$ Renormalization  $\alpha_m:=x^{-1}\tilde{\alpha}_m$ gives the desired system.
\end{proof}

\begin{remark}\label{rem:fundamentalsystem}
 Since $H$ is also a closed subgroup of the absolute Galois group $G_L$ of $L$ it possesses   a countable  fundamental system $(H_m)_m$ of open neighbourhoods   of the identity, as for any $n>0$ the local field $L$   of characteristic $0$ has only finitely many extensions of degree smaller than $n$.
\end{remark}

\begin{proof}
The latter statement reduces easily to finite Galois extensions $L'$ of $L$, which are known to be solvable, i.e. $L'$ has a series of at most $n$ intermediate fields $L\subseteq L_1\subseteq \ldots \subseteq L_n=L'$ such that each subextension is abelian. Now its known by class field theory that each local field in characteristic $0$ only has finitely many abelian extensions of a given degree.
\end{proof}

We write $\mathcal{C}^n(G,V)$ for the abelian group of continuous $n$-cochains of a profinite group $G$ with values in a topological abelian group $V$ carrying a continuous $G$-action and $\partial$ for the usual differentials. In particular, we endow $\mathcal{C}^n(H,W)$ with the maximum norm $\|-\|$ and consider the subspace $\mathcal{C}^n(H,W)^\delta:=\bigcup_{H'\unlhd H \mbox{ open }}\mathcal{C}^n(H/H',W)\subseteq \mathcal{C}^n(H,W)$ of those cochains with are even continuous with respect to the discrete topology of $W.$

\begin{lemma}\phantomsection\label{lem:homotopy}
\begin{enumerate}
\item The completion of ${\mathcal{C}^n(H,W)^\delta}$ with respect to the maximum norm equals $\mathcal{C}^n(H,W).$
    %More precisely, for every $f\in {\mathcal{C}^n(H,W)}$ there exists a sequence of open normal subgroups $H_m$ of $H$  such that $f$ belongs to the completion of $\mathcal{C}^n(H,(H_m),W):=\bigcup_{m}\mathcal{C}^n(H/H_m,W)$.
\item There exist $(\mathbb{C}_p^\flat)^H$-linear continuous maps
\[\sigma^n:\mathcal{C}^n(H,W)\to \mathcal{C}^{n-1}(H,W) \] satisfying $\|f-\partial\sigma^n f\|\leq c\|\partial f\|.$
\end{enumerate}
\end{lemma}

\begin{proof}
 Since the space $\mathcal{C}^n(H,W) $ is already complete we only have to show that an arbitrary cochain $f$ in it can be approximated by a Cauchy sequence $f_m$  in $\mathcal{C}^n(H,W)^\delta$. To this end we observe that, given any $m$, the induced cochain $H^n\xrightarrow{f}W\xrightarrow{pr_m} W/W_{\geq m}$ comes, for some open normal subgroup $H_m$, from a cochain in $\mathcal{C}^n(H/H_m,W/W_{\geq m}),$ which in turn gives rise to $f_m \in \mathcal{C}^n(H,W)^\delta$ when composing with any set theoretical section $W/W_{\geq m}\xrightarrow{s_m}W$ of the canonical projection $W\xrightarrow{pr_m} W/W_{\geq m}.$ Note that $s_m$ is automatically continuous, since $W/W_{\geq m}$ is discrete.  By construction we have $\|f-f_m\|\leq \frac{1}{p^m}$ and $(f_m)_m$ obviously is a Cauchy sequence. This shows (i).

 For (ii)  recall from Lemma \ref{lem:TS1} together with Remark \ref{rem:fundamentalsystem} the existence  of a trace compatible system $(y_{H'})_{H'}$ of elements $y_{H'}\in (\mathbb{C}_p^\flat)^{H'}$ with $|y_{H'}|_\flat < c$  and $Tr_{H|H'}(y_{H'} )=1$, where $H'$ runs over the open normal subgroups of $H$.  Now we first define  $(\mathbb{C}_p^\flat)^H$-linear  maps
\[
\sigma^n: \mathcal{C}^n(H,W)^\delta \to \mathcal{C}^{n-1}(H,W)
\]
satisfying $\|f-\partial\sigma^n f\|\leq c\|\partial f\|$ and $\| \sigma^n f\|\leq c \|f\|$ by setting for $f\in  \mathcal{C}^n(H/H',W)$
\[
\sigma^n(f):=y_{H'}\cup f
\]
(by considering $y_{H'}$ as a $-1$-cochain), i.e.,
\[
\sigma^n(f)(h_1,\ldots, h_{n-1})=(-1)^n\sum_{\tau\in H/H'} (h_1 \ldots h_{n-1}\tau)(y_{H'}) f(h_1,\ldots, h_{n-1},\tau).
\]
The inequality $\|y_{H'}\cup f\|\leq c\|f\|$ follows immediately from this description, see the proof of \cite[Lem.\ 14.3.1.]{BC}.
Upon noting that $\partial y_{H'}=Tr_{H|H'}(y_{H'})=1$, the Leibniz rule for the differential $\partial$ with respect to the cup-product then implies that
\[
f-\partial(y_{H'}\cup f)=y_{H'}\cup \partial f,
\]
hence
\[
\|f-\partial(y_{H'}\cup f)\| \leq c \|\partial f\|
\]
by the previous inequality, see again loc.\ cit. In order to check that this map $\sigma^n$ is well defined we assume that $f$ arises also from a cochain in  $\mathcal{C}^n(H/H'',W)$. Since we may make the comparison within $\mathcal{C}^n(H/(H'\cap H''),W)$ we can assume without loss of generality that $H''\subseteq H'$. Then
\begin{align*}
  (y_{H''}\cup f)( h_1,\ldots, h_{n-1})&=(-1)^n\sum_{\tau\in H/H''} (h_1 \ldots h_{n-1}\tau)(y_{H''}) f(h_1,\ldots, h_{n-1},\tau) \\
    & =(-1)^n\sum_{\tau\in H/H'} \left(h_1 \ldots h_{n-1}\sum_{\tau'\in H'/H''}\tau'\right)(y_{H''}) f(h_1,\ldots, h_{n-1},\tau)\\
    & =(-1)^n\sum_{\tau\in H/H'} \left(h_1 \ldots h_{n-1}\right)(\sum_{\tau'\in H'/H''}\tau'(y_{H''})) f(h_1,\ldots, h_{n-1},\tau)\\
        & =(-1)^n\sum_{\tau\in H/H'} \left(h_1 \ldots h_{n-1}\right)( y_{H'}) f(h_1,\ldots, h_{n-1},\tau)\\
        &= (y_{H'}\cup f)( h_1,\ldots, h_{n-1})
\end{align*} using the trace compatibility in the fourth equality. Finally the inequality $\| \sigma^n f\|\leq c \|f\|$ implies that $\sigma^n$ is continuous on $\mathcal{C}^n(H,W)^\delta$ and therefore extends continuously to its completion $\mathcal{C}^n(H,W)$.
\end{proof}

The proof of  Prop.\  \ref{prop:TatevanishingH1} is now an immediate consequence of   Lemma  \ref{lem:homotopy}(ii).

\Footnote{\begin{remark} \label{rem:proetale} Apparently, the perfectoid theory for pro\'{e}tale sites gives vanishing statements for higher cohomology like in the last Prop.\  almost for free: Indeed, let $X$ be the adic space $ Spa(L,o_L) $. Then $Y_\infty:=Spa( \widehat{L_\infty},o_{\widehat{L_\infty}})$ forms a perfectoid subdomain of $X_{pro\acute{e}t}$, whence \cite[Lem.\ 9.3.4]{KLI} claims that any $\mathcal{P}$ of the period sheaves of Definition 9.3.3 in (loc.\ cit.) is acyclic on $Y_\infty,$ i.e. $H^i((Y_\infty)_{pro\acute{e}t},\mathcal{P})=0$ for $i>0.$ Using the (degenerating) Cartan-Leray spectral sequence       (as in the proof of \cite[Prop.\ 3.7 (iii), Prop.\ 6.16 (i) - see  also Thm.\ 4.9, Lem.\ 4.10, Cor.\ 6.6]{Scho}- unfortunately the identification of Cech-cochains of the covering below with continuous maps from $H_L^n$ to the required coefficients requires further vanishing results and this approach is not as straight forward as we hoped for!)
\[H^p(H_L,H^0(Z_{pro\acute{e}t},\mathcal{P}))= H^{p}((Y_\infty)_{pro\acute{e}t},\mathcal{P})\]
for the for the universal cover, i.e., pro-(finite)\'{e}tale  map  $Z:=Spa(\mathbb{C}_p, o_{\mathbb{C}_p})\to Y_\infty$ mit Galois group $H_L$ (as the categories of finite \'{e}tale objects over $Spec({L_\infty} )$, $Spec(\widehat{L_\infty} )$ and $Spa( \widehat{L_\infty},o_{\widehat{L_\infty}})$ are equivalent) and using that
\begin{equation*}
  H^q(Z\times_{Y_\infty}\cdots\times_{Y_\infty} Z ,\mathcal{P})=                                                                 0\hbox{ for } q>0
\end{equation*}
we obtain that $H^i(H_L,\mathcal{P}(Z))=0$ for $i>0.$ In particular, this holds for the perfect period rings $\mathcal{P}(Z)\in\{\tilde{\cR}, \tilde{\mathbf{A}}, \tilde{\mathbf{A}}^\dagger,  \tilde{\mathbf{B}}^\dagger \}.$
  The case $\mathbb{C}_p^\flat $ corresponds via the perfectoid correspondence to the aclycity of the completed structure sheaf \cite[Lem.\ 9.2.8]{KLI}.
\end{remark}}
\Footnote{$H^1(H_L,\tilde{\cR})=0$  seems to be equivalent to $\tilde{\cR}_L\subseteq\tilde{\cR}$ being faithfully flat. We do not need this statement, if we use Lemma \ref{lem:stronghypKed} instead! How can one prove these results avoiding the pro\'{e}tale site?}

\section{The  functors \texorpdfstring{$D$}{D},  \texorpdfstring{$\tilde{D}$}{tilde D} and \texorpdfstring{$\tilde{D}^\dagger$}{tilde D dagger}}\label{subsubsec:D}

     Let $\Rep_{o_L}(G_L),$ $\Rep_{o_L,f}(G_L)$ and $\Rep_{L}(G_L) $ denote the category of finitely generated  $o_L$-modules, finitely generated  free $o_L$-modules and finite dimensional $L$-vector spaces, respectively, equipped with a continuous linear $G_L$-action.
The following result is established in \cite[ Thm.\ 1.6]{KR} (see also \cite[Thm.\ 3.3.10]{GAL}) and \cite[Prop.\ 4.4 (ii)]{SV15}.

\begin{theorem}\label{KR-equiv}
The functors
\begin{equation*}
    T \longmapsto D(T) := (\mathbf{A} \otimes_{o_L} T)^{H_{L}} \qquad\text{and}\qquad M \longmapsto (\mathbf{A} \otimes_{\mathbf{A}_L} M)^{\varphi_q \otimes \varphi_M =1}
\end{equation*}
are exact quasi-inverse equivalences of categories between $\Rep_{o_L}(G_L)$ and the category $\mathfrak{M}^{et}(\mathbf{A}_L)$ of finitely generated \'{e}tale $\varphi_L,\Gamma_L)$-modules over $ \mathbf{A}_L.$ Moreover, for any $T$ in $\Rep_{o_L}(G_L)$ the natural map \begin{equation}\label{f:compD}
 \mathbf{A} \otimes_{\mathbf{A}_L} D(T) \xrightarrow{\; \cong \;} \mathbf{A} \otimes_{o_L} T
\end{equation}
is  an isomorphism (compatible with the $G_L$-action and the Frobenius on both sides).
\end{theorem}

In the following we would like to establish a   version of the above for $\tilde{\mathbf{A}}$ and prove similar properties for it.  In the classical situation such versions have been studied by Kedlaya et al using the unramified rings of Witt vectors $W(R)$. In our Lubin-Tate situation we have to work with ramified Witt vectors $W(R)_L$. Many results and their proofs transfer almost literally from the classical setting. Often we will try to at least sketch the proofs for the convenience of the reader, but when we just quote results from the classical situation, e.g. from \cite{KLI}, this usually means that the transfer is purely formal.

We start defining  $\tilde{\mathbf{A}}:=W(\mathbb{C}_p^\flat)_L$ and
\[
\tilde{\mathbf{A}}^\dagger:=\{x=\sum_{n\geq 0}\pi_L^n[x_n]\in\tilde{\mathbf{A}}: |\pi_L^{n}\|x_n|_\flat^r \xrightarrow{n \rightarrow \infty} 0\mbox{ for some } r>0\}
\]
as well as $\tilde{D}(T):=(\tilde{\mathbf{A}}\otimes_{o_L} T)^{H_L}$ and $\tilde{D}^\dagger(T):=(\tilde{\mathbf{A}}^\dagger\otimes_{o_L} T)^{H_L}.$\Footnote{
In the literature one also finds the subring $ \tilde{\mathbf{A}}^\dagger_{\leq 1}:=\bigcup_{r>0}W^r_{\leq 1}(\mathbb{C}_p^\flat)_L$ of $ \tilde{\mathbf{A}}^\dagger $ where $W^r_{\leq 1}(\mathbb{C}_p^\flat)_L=\{ x\in W^r(\mathbb{C}_p^\flat)_L|\; |x|_{r}\leq 1 \}$ consists of those $x\in \tilde{\mathbf{A}}$ such that   $|\pi_L^{n}\|x_n|_\flat^r \xrightarrow{n \rightarrow \infty} 0$ and $|\pi_L^{n}\|x_n|_\flat^r\leq 1$ for all $n.$ Denoting by $\tilde{\mathbf{A}}^{\dagger,s}_{St}$ the ring defined in \cite[Def.\ 3.4]{Ste} we have the equality $W^r_{\leq 1}(\mathbb{C}_p^\flat)_L = \tilde{\bf A}^{\dagger,\frac{q-1}{qr}}_{St}$ corresponding to $\tilde{\bf A}^{\dagger}_{(\frac{q-1}{qr})^-}$ in the notation of  \cite[\S II.1]{ChCo} for $q=p$. For these relations use the following  normalisations compatible with \cite{Ste}:
 $|\pi_L|=\frac{1}{q}$, $v_{\mathbb{C}_p^\flat}(\omega)=\frac{q}{q-1},$ $v_{\pi_L} (\pi_{L})=1,$ $|x|_\flat=q^{-v_{\mathbb{C}_p^\flat}}, $ $|\omega|_\flat=q^{-\frac{q}{q-1}}=|\pi_L|^{\frac{q}{q-1}},$  where $\omega=\omega_{LT}\mod \pi_L$.
Furthermore, $|x|_{r}=|\pi_L|^{V(x,\frac{q-1}{rq})}$ and $|\omega_{LT}|_{r}=|\pi_L|^{\frac{rg}{q-1}}=|{\omega}|^r_\flat,$ where $V(x,r):=\inf_k\left(v_{\mathbb{C}_p^\flat}(x_k)\frac{q-1}{rq}+k\right)$ for $x=\sum_{k\geq 0} \pi_L^k[x_k]\in \tilde{\bf A}.$ For $x\in \tilde{\mathbf{A}}^\dagger$ we have $V(x,r)=\frac{q-1}{rq}V_{St}(x,r)$ where $V_{St}(x,r) $ uses the notation in \cite{Ste}.
Note also that $\omega_{LT}^{-1}$ is contained in $W^{\frac{q-1}{q}}_{\leq 1}(\widehat{L_\infty}^\flat)_L$ by  \cite[Lem.\ 3.10]{Ste} (in analogy with \cite[Cor. II.1.5]{ChCo}).

}

More generally, let $K$ be   any perfectoid field containing $L$ and let $K^\flat$ denote its tilt. For $r>0$ let $W^r(K^\flat)_L$ be the set of $x=\sum_{n=0}^{\infty} \pi_L^n[x_n]\in W(K^\flat)_L$ such that $|\pi_L|^n|x_n|^r_\flat$ tends to zero as $n$ goes to $\infty.$ This is a subring by \cite[Prop.\ 5.1.2]{KLI} on which the function
\[
|x|_{r}:=\sup_n \{|\pi_L^{n}\| x_n|^r_\flat\}=\sup_n \{q^{-n}| x_n|^r_\flat\}
\]
is a complete multiplicative norm; it extends multiplicatively to $W^r(K^\flat)_L[\frac{1}{\pi_L}].$
Furthermore, $W^\dagger(K^\flat)_L:=\bigcup_{r>0}W^r(K^\flat)_L$ \footnote{In \cite{Ked05} it is denoted by $W^\dagger(K^\flat)_L$.  } is a  henselian discrete valuation ring by \cite[Lem.\ 2.1.12]{Ked05}, whose $\pi_L$-adic completion equals $W(K^\flat)_L$ since they coincide modulo $\pi_L^n$. Then $\tilde{\bf A}^{\dagger}=W^\dagger(\mathbb{C}_p^\flat)_L$, and we write $\tilde{\bf A}_L$ and $\tilde{\bf A}_L^{\dagger}$   for $W(\hat{L}_\infty^\flat)_L$ and $W^\dagger(\hat{L}_\infty^\flat)_L$, respectively.  We set $\tilde{\bf B}_L=\tilde{\bf A}_L[\frac{1}{\pi_L}]$,  $\tilde{\bf B}=\tilde{\bf A}[\frac{1}{\pi_L}]$, $\tilde{\bf B}_L^{\dagger}=\tilde{\bf A}_L^{\dagger}[\frac{1}{\pi_L}]$ and $\tilde{\bf B}^{\dagger}=\tilde{\bf A}^{\dagger}[\frac{1}{\pi_L}]$ for the corresponding fields of fractions.

\begin{remark}
By the Ax-Tate-Sen theorem \cite{Ax} and since $\mathbb{C}_p^\flat$ is the completion of an algebraic closure $\overline{\hat{L}_\infty^\flat}$ he have that
$(\mathbb{C}_p^\flat)^H=((\overline{\hat{L}_\infty^\flat})^H)^\wedge$
for any closed subgroup $H\subseteq H_L$, in particular $(\mathbb{C}_p^\flat)^{H_L} = \hat{L}_\infty^\flat$.  As completion of an algebraic extension   of the perfect field $\hat{L}_\infty^\flat$ the field $(\mathbb{C}_p^\flat)^H$ is perfect, too. Moreover, we have $\tilde{\mathbf{A}}^{H_L}=\tilde{\mathbf{A}}_L$,  $(\tilde{\mathbf{A}}^\dagger)^{H_L}=\tilde{\mathbf{A}}_L^\dagger$ and analogously for the rings $\tilde{\bf B}$ and $\tilde{\bf B}^{\dagger}$.  It also follows that $\tilde{\mathbf{A}}$ is the $\pi_L$-adic completion of a maximal unramified extension of $\tilde{\mathbf{A}}_L$.
\end{remark}

\Footnote{ We expect that \cite[Lem.\ 3.18, Prop.\ 3.7 (iii)]{Scho} implies that  $H^n(H,\tilde{\mathbf{A}})=0$ for all $n\geq 1$ and $H\subseteq H_L$  any closed subgroup. }

\begin{lemma}
The rings $\mathbf{A}_L$ and $\mathbf{A}$ embed into $\tilde{\mathbf{A}}_L$ and $\tilde{\mathbf{A}}$, respectively.
\end{lemma}

\begin{proof}
 The embedding $\mathbf{A}_L \hookrightarrow \tilde{\mathbf{A}}_L$ is explained  in \cite[p.\ 94]{GAL}.
%\footnote{
%Note that $\omega_\varphi=\{[\iota(t)]\}\in \tilde{\bf A}_\psi^+:=W(o_{\hat{L}_\infty}^\flat)_L\subseteq\tilde{\bf A}_L$   (in the notation of \cite{GAL}).
%Since  $\omega_\varphi\equiv [\bar{\omega}]\mod \pi_L$, $\omega_\varphi$ is invertible in $\tilde{\bf A}_L $ and it is clear that $o_L[[\omega_\varphi]][\frac{1}{\omega_\varphi}]\subseteq  \tilde{\mathbf {A}}, $ whence $\mathbf{A}_L\subseteq \tilde{\mathbf{A}}_L.$}
Moreover, $\mathbf{A}$ is the $\pi_L$-adic completion of the maximal unramified extension of $\mathbf{A}_L$ inside $\tilde{\mathbf{A}} = W(\mathbb{C}_p^\flat)_L$ (cf.\ \cite[\S3.1]{GAL}).
\end{proof}

On $\tilde{\mathbf{A}}=W(\mathbb{C}_p^\flat)_L$ the weak topology is defined to be the product topology of the valuation topologies on the components $\mathbb{C}_p^\flat.$ The induced topology on any subring $R$ of it is also called weak topology of $R$. If $M$ is a finitely generated $R$-module, then we call the {\it canonical} topology of $M$  {\it (with respect to the weak topology of $R$)} the quotient topology with respect to any surjection $R^n\twoheadrightarrow M$ where the free module carries the product topology; this is independent of any  choices.
We recall that a $(\varphi_L,\Gamma_L)$-module $M$ over $R\in\{{\mathbf{A}}_L, \tilde{\mathbf{A}}_L, \tilde{\mathbf{A}}_L^\dagger\} $ is a finitely generated $R$-module $M$ together with
\begin{itemize}
  \item[--] a $\Gamma_L$-action on $M$ by semilinear automorphisms which is continuous for the weak topology
%\footnote{In case $L = \mathbb{Q}_p$ we have automatic continuity. The simplest instance of this is the fact that any abstract group homomorphism $\Gamma_{\mathbb{Q}_p} \longrightarrow \mathbb{Z}_p^\times$ is continuous. By restricting to sufficiently small open subgroups this reduces to the claim that any abstract group homomorphism $\mathbb{Z}_p \longrightarrow \mathbb{Z}_p$ is continuous, i.e., is determined by its value in $1$. This follows from the triviality of any group homomorphism $\mathbb{Z}_p/\mathbb{Z} \longrightarrow \mathbb{Z}_p$. The latter holds because, by the surjectivity of the projection map $\mathbb{Z} \longrightarrow \mathbb{Z}_p/p \mathbb{Z}_p$, the group $\mathbb{Z}_p/\mathbb{Z}$ is $p$-divisible. Actually we always have automatic continuity, see \cite{GAL} Thm.\ 2.2.8.},
and
  \item[--] a $\varphi_L$-linear endomorphism $\varphi_M$ of $M$ which commutes with the $\Gamma_L$-action.
\end{itemize}
We let $\mathfrak{M}(R)$ denote the category of   $(\varphi_L,\Gamma_L)$-modules $M$ over $R$. Such a module $M$  is called \'{e}tale if the linearized map
\begin{align*}
     \varphi_M^{lin} : R \otimes_{R,\varphi_L} M & \xrightarrow{\; \cong \;} M \\
                                           f \otimes m & \longmapsto f \varphi_M (m)
\end{align*}
is bijective. We let   $\mathfrak{M}^{\acute{e}t}(R)$  denote  the full subcategory of \'{e}tale $(\varphi_L,\Gamma_L)$-modules over $R$.

\begin{definition}\label{def:etale}
 For $*={\mathbf{B}}_L, \tilde{\mathbf{B}}_L, \tilde{\mathbf{B}}_L^\dagger$   we write $ \mathfrak{M}^{\acute{e}t}(*):= \mathfrak{M}^{\acute{e}t}(*')\otimes_{o_L}L$ with $*'={\mathbf{A}}_L, \tilde{\mathbf{A}}_L, \tilde{\mathbf{A}}_L^\dagger$, respectively, and call the objects \'{e}tale $(\varphi_L,\Gamma_L)$-modules over $*.$
\end{definition}

\begin{lemma}\label{lem:contGammaAction}
Let $G$ be a profinite group and $R\to S$ be a topological monomorphism of topological $o_L$-algebras, for which there exists a system of open neighbourhoods of $0$ consisting of $o_L$-submodules. Consider   a finitely generated $R$-module $M$, for which  the canonical map $M\to S\otimes_R M$ is injective (e.g.\ if $S$ is faithfully flat over $R$ or $M$ is free, in addition), and  endow it  with the canonical topology with respect to $R.$ Assume that $G$ acts continuously, $o_L$-linearly and compatible on $R$ and $S$ as well as continuously and $R$-semilinearly on $M$.  Then the diagonal $G$-action on $S\otimes_R M$   is continuous with regard to the canonical topology with respect to $S.$
\end{lemma}

\begin{proof}
Imitate the proof of \cite[Lem.\ 3.1.11]{GAL}.
\end{proof}

\begin{proposition}\label{Dbasechange}
 The canonical map
 \begin{equation}\label{f:Dbasechange}
    \tilde{\mathbf{A}}_L\otimes_{ {\mathbf{A}}_L} {D}(T)\xrightarrow{\cong} \tilde{D}(T)
 \end{equation} is an isomorphism and
 the functor   $\tilde{D}(-):\Rep_{o_L}(G_L) \to \mathfrak{M}^{\acute{e}t}(\tilde{\mathbf{A}}_L)$ is exact. Moreover, we have a comparison isomorphism
     \begin{equation}\label{f:comptildeD}
      \tilde{\mathbf{A}}\otimes_{ \tilde{\mathbf{A}}_L} \tilde{D}(T) \xrightarrow{\cong}  \tilde{\mathbf{A}}\otimes_{o_L} T.
     \end{equation}
\end{proposition}

\begin{proof}
The isomorphism \eqref{f:Dbasechange} implies formally the isomorphism \eqref{f:comptildeD} after base change of the comparison isomorphism \eqref{f:compD}.  Secondly, the isomorphism \eqref{f:Dbasechange}, resp.\ \eqref{f:comptildeD}, implies easily that $\tilde{D}(T)$ is finitely generated, resp.\ \'{e}tale. Thirdly, since the ring extension $\tilde{\mathbf{A}}_L/ { {\mathbf{A}}_L}$ is  faithfully flat as local extension of (discrete) valuation rings, the exactness of $\tilde{D}$ follows from that of $D.$  Moreover,  the isomorphism \eqref{f:Dbasechange} implies by Lemma \ref{lem:contGammaAction} that $\Gamma_L$ acts continuously on $\tilde{D}(T)$, i.e., the functor $\tilde{D}$ is well-defined.  Thus we only have to prove that
\begin{equation*}
    \tilde{\mathbf{A}}_L\otimes_{ {\mathbf{A}}_L}  (\mathbf{A}\otimes_{ o_L}T)^{H_L}\xrightarrow{\cong}  (\tilde{\mathbf{A}}\otimes_{ o_L}T)^{H_L}
 \end{equation*}
s an isomorphism. To this aim let us {\it assume first that $T$ is finite.} Then we find an open normal subgroup $H\unlhd H_L$ which acts trivially on $T.$ Application of the subsequent Lemma \ref{lem:Ginv} to $M=(\mathbf{A}\otimes_{ o_L}T)^{H} $ and $G=H_L/H$ interprets the left hand side as $\left( \tilde{\mathbf{A}}_L\otimes_{ {\mathbf{A}}_L}(\mathbf{A}\otimes_{ o_L}T)^{H}  \right)^{H_L/H}$ while the right hand side equals $\left(  (\tilde{\mathbf{A}}\otimes_{ o_L}T)^{H}  \right)^{H_L/H}$. Hence it suffices to establish the isomorphism
\begin{equation*}
    \tilde{\mathbf{A}}_L\otimes_{ {\mathbf{A}}_L}  (\mathbf{A}\otimes_{ o_L}T)^{H}\xrightarrow{\cong}  (\tilde{\mathbf{A}}\otimes_{ o_L}T)^H.
 \end{equation*}
By Lemma \ref{lem:Hinv} below this is reduced to showing that the canonical map
\begin{equation*}
    \tilde{\mathbf{A}}_L\otimes_{ {\mathbf{A}}_L}  \mathbf{A}^{H}\otimes_{ o_L}T\xrightarrow{\cong}  \tilde{\mathbf{A}}^H\otimes_{ o_L}T
 \end{equation*}
is an isomorphism, which follows from   Lemma \ref{lem:lineardisjoint} below.
{\it Finally let $T$ be arbitrary.} Then we have isomorphisms
\begin{align*}
 \tilde{\mathbf{A}}_L\otimes_{ {\mathbf{A}}_L} {D}(T)   & \cong \tilde{\mathbf{A}}_L\otimes_{ {\mathbf{A}}_L} \varprojlim_{n}{D}(T/\pi_L^nT) \\
    & \cong \tilde{\mathbf{A}}_L\otimes_{ {\mathbf{A}}_L} \varprojlim_{n}{D}(T)/\pi_L^n{D}(T) \\
    & \cong   \varprojlim_{n}\tilde{\mathbf{A}}_L\otimes_{ {\mathbf{A}}_L}{D}(T)/\pi_L^n{D}(T)  \\
    & \cong  \varprojlim_{n}\tilde{\mathbf{A}}_L\otimes_{ {\mathbf{A}}_L}{D}(T/\pi_L^n T)   \\
    & \cong  \varprojlim_{n} \tilde{D}(T/\pi_L^n T) \\
    &\cong \tilde{D}(T),
\end{align*}
where we use for the second and fourth equation exactness of $D$, for the second last one the case of finite $T$ and for the first, third and last equation the elementary divisor theory for the discrete valuation rings $o_L$, $\mathbf{A}_L$ and $\tilde{\mathbf{A}}_L,$ respectively.
\end{proof}

\begin{lemma}\label{lem:Ginv}
Let $A\to B$ be a flat extension of rings and $M$ an $A$-module with an $A$-linear action by a finite group $G$. Then $B\otimes_A M$ carries a $B$-linear $G$-action and we have
\[\left(B\otimes_A M\right)^G=B\otimes_A M^G.\]
\end{lemma}

\begin{proof}
Apply the exact functor $B\otimes_A -$ to the exact sequence
\[\xymatrix@C=0.5cm{
  0 \ar[r] & M^G \ar[rr]^{ } && M \ar[rr]^-{(g-1)_{g\in G}} && \bigoplus_{g\in G} M  ,}\] which gives the desired description of $\left(B\otimes_A M\right)^G.$
\end{proof}

\begin{lemma}\label{lem:Hinv}
Let $A$ be $\mathbf{A}$,  $\mathbf{A}_L^{nr}$, $\tilde{\mathbf{A}}^\dagger$ or $\tilde{\mathbf{A}}$ and $T$ be a finitely generated $o_L$-module with trivial action by an open subgroup $H\subseteq H_L$. Then $(A\otimes_{o_L}T)^H=A^H\otimes_{o_L} T.$ Moreover, $\mathbf{A}^{H} $ and $\tilde{\mathbf{A}}^{H} $ are   free $\mathbf{A}_{L}$- and $\tilde{\mathbf{A}}_{L}$-modules of finite rank, respectively.
\end{lemma}

\begin{proof}
  Since $T\cong \bigoplus_{i=1}^r o_L/\pi_L^{n_i}o_L$ with $n_i \in \mathbb{N} \cup \{\infty\}$ we may assume that $T=o_L/\pi_L^no_L$ for some $n \in \mathbb{N} \cup \{\infty\}$. We then we have to show that
\begin{align}\label{f:Hinmodpi}
     (A/\pi_L^n A)^H= &  A^H/\pi_L^nA^H
\end{align}
For $n = \infty$ there is nothing to prove.

\textit{The case $n = 1$:}
First of all we have $\mathbf{A}/\pi_L\mathbf{A} =\mathbf{A}_L^{nr}/\pi_L\mathbf{A}_L^{nr}=\mathbf{E}_L^{sep}$. On the other hand, by the Galois correspondence between unramified extensions and their residue extensions, we have that $(\mathbf{E}_L^{sep})^H$ is the residue field of $(\mathbf{A}_L^{nr})^{H}$. Hence the case $n=1$ holds true for $A = \mathbf{A}_L^{nr}$.  After having finished all cases for $A = \mathbf{A}_L^{nr}$ we will see at the end of the proof that $(\mathbf{A}_L^{nr})^{H}  =\mathbf{A}^{H}$. Therefore the case $n=1$ for $A = \mathbf{A}$ will be settled, too.

For $A = \tilde{\mathbf{A}}$ we only need to observe that $\tilde{\mathbf{A}}/\pi_L\tilde{\mathbf{A} } =W(\mathbb{C}_p^\flat)_L/\pi_L W(\mathbb{C}_p^\flat)_L=\mathbb{C}_p^\flat$ and that $(\mathbb{C}_p^\flat)^H$ is the residue field of $(W(\mathbb{C}_p^\flat)_L)^{H}=W((\mathbb{C}_p^\flat)^{H})_L.$

For $A = \tilde{\mathbf{A}}^\dagger$ we argue by the following commutative diagram
  \[\xymatrix{
    (\mathbf{C}_p^\flat)^H \ar[dr]_-{\cong } \ar[r]^-{\cong } & W^\dagger((\mathbf{C}_p^\flat)^H )_L / \pi_L W^\dagger((\mathbf{C}_p^\flat)^H )_L \ar[rr]^-{\cong } & &(\tilde{\mathbf{A}}^\dagger)^H/\pi_L (\tilde{\mathbf{A}}^\dagger)^H \ar[d]^{ } \\
        & \tilde{\mathbf{A}}^H/\pi_L \tilde{\mathbf{A}}^H \ar[r]^-{\cong } & (\tilde{\mathbf{A}}/\pi_L \tilde{\mathbf{A}})^H  \ar[r]^-{\cong } &(\tilde{\mathbf{A}}^\dagger/\pi_L \tilde{\mathbf{A}}^\dagger)^H .  }\]

\textit{The case $1 < n < \infty$:} This follows by induction using the commutative diagram with exact lines
\[\xymatrix{
    0  \ar[r]^{ } & A^H/\pi_L^nA^H \ar[d]_{\cong } \ar[r]^{\pi_L \cdot} & A^H/\pi_L^{n+1}A^H \ar[d]_{ } \ar[r]^{ } &  A^H/\pi_LA^H \ar[d]_{\cong} \ar[r]^{ } & 0 \\
    0 \ar[r]^{ } & (A/\pi_L^n A)^H \ar[r]^{\pi_L \cdot} & (A/\pi_L^{n+1} A)^H \ar[r]^{ } & (A/\pi_L A)^H ,}
\]
in which the outer vertical arrows are isomorphism by the case $n=1$ and the induction hypothesis.

Finally we can check, using the above equality \eqref{f:Hinmodpi} for $A = \mathbf{A}_L^{nr}$ in the third equation:
\begin{align*}
    \mathbf{A}^{H} & = \left( \varprojlim_{n} \mathbf{A}_L^{nr}/\pi_L^n\mathbf{A}_L^{nr}\right)^H\\
     & = \varprojlim_{n} \left(\mathbf{A}_L^{nr}/\pi_L^n\mathbf{A}_L^{nr}\right)^H \\
      & = \varprojlim_{n} \left(\mathbf{A}_L^{nr})^H/\pi_L^n(\mathbf{A}_L^{nr}\right)^H  \\
      &  = (\mathbf{A}_L^{nr})^H.
\end{align*}
Note that $(\mathbf{A}_L^{nr})^H$ is a finite unramified extension of $\mathbf{A}_L$ and therefore is $\pi_L$-adically complete. We also see that $\mathbf{A}^{H} $ is a free $\mathbf{A}_{L}$-module of finite rank. Similarly,  $W(\mathbb{C}_p^\flat)_L^{H}  \cong (W(\hat{L}_\infty^\flat)_L^{nr})^H$ is a free $W(\hat{L}_\infty^\flat)_L$-module of finite rank.
\end{proof}

\begin{lemma}\label{lem:lineardisjoint}
 For any open subgroup $H$ of $H_L$ the canonical maps \begin{align*}
    W(\hat{L}_\infty^\flat)_L\otimes_{ {\mathbf{A}}_L}  \mathbf{A}^{H} &\xrightarrow{\cong}   W((\mathbb{C}_p^\flat)^H)_L, \\
    W(\hat{L}_\infty^\flat)_L\otimes_{ \tilde{\mathbf{A}}_L^\dagger}  (\tilde{\mathbf{A}}^\dagger)^{H} &\xrightarrow{\cong}  W((\mathbb{C}_p^\flat)^H)_L
 \end{align*}
 are  isomorphisms.
\end{lemma}

\begin{proof}
We begin with the first isomorphism. Since $\mathbf{A}^{H} $ is finitely generated free over $\mathbf{A}_{L}$ by Lemma \ref{lem:Hinv}, we have
\begin{equation*}
  W(\hat{L}_\infty^\flat)_L\otimes_{ {\mathbf{A}}_L}  \mathbf{A}^{H}  \cong \left(\varprojlim_{n} W_n(\hat{L}_\infty^\flat)_L\right)\otimes_{ {\mathbf{A}}_L}  \mathbf{A}^{H}
     \cong  \varprojlim_{n} \left(W_n(\hat{L}_\infty^\flat)_L\otimes_{ {\mathbf{A}}_L}  \mathbf{A}^{H}\right) .
\end{equation*}
It therefore suffices to show the corresponding assertion for Witt vectors of finite length:
\begin{equation*}
    W_n(\hat{L}_\infty^\flat)_L\otimes_{ {\mathbf{A}}_L}  \mathbf{A}^{H}/ \pi_L^n \mathbf{A}^H = W_n(\hat{L}_\infty^\flat)_L\otimes_{ {\mathbf{A}}_L}  \mathbf{A}^{H} \xrightarrow{\cong}  W_n((\mathbb{C}_p^\flat)^H)_L .
 \end{equation*}
  To this aim we first consider the case $n=1$. From \eqref{f:Hinmodpi} we know that $\mathbf{A}^{H}/ \pi_L^n \mathbf{A}^H = (\mathbf{E}_L^{sep})^H$. Hence we need to check that
\begin{equation*}
  \hat{L}_\infty^\flat \otimes_{\mathbf{E}_L} (\mathbf{E}_L^{sep})^H \xrightarrow{\cong} (\mathbb{C}_p^\flat)^H
\end{equation*}
is an isomorphism. Since the perfect hull $\mathbf{E}_L^{perf}$ of $\mathbf{E}_L$ (being purely inseparable and normal) and $(\mathbf{E}_L^{sep})^H$ (being separable) are linear disjoint extensions of $\mathbf{E}_L$ their tensor product is equal to the composite of fields $\mathbf{E}_L^{perf} (\mathbf{E}_L^{sep})^H$ (cf.\ \cite[Thm.\ 5.5, p.\ 188]{Coh}), which moreover has to have degree $[H_L : H]$ over $\mathbf{E}_L^{perf}$. Since the completion of the tensor product is $\hat{L}_\infty^\flat \otimes_{\mathbf{E}_L} (\mathbf{E}_L^{sep})^H$, we see that the completion of the field $\mathbf{E}_L^{perf} (\mathbf{E}_L^{sep})^H$ is the composite of fields $\hat{L}_\infty^\flat (\mathbf{E}_L^{sep})^H$, which has degree $[H_L : H]$ over $\hat{L}_\infty^\flat$. But $\hat{L}_\infty^\flat (\mathbf{E}_L^{sep})^H \subseteq (\mathbb{C}_p^\flat)^H$. By the Ax-Tate-Sen theorem $(\mathbb{C}_p^\flat)^H$ has also degree $[H_L : H]$ over $\hat{L}_\infty^\flat$. Hence the two fields coincide, which establishes the case $n=1$.

The commutative diagram
\[ \xymatrix{
   \hat{L}_\infty^\flat \otimes_{ {\mathbf{A}}_L}  \mathbf{A}^{H} \ar[d]_{\varphi_q^m \otimes \id}^{\cong} \ar[r]^-{\cong} & (\mathbb{C}_p^\flat)^H \ar[d]^{\varphi_q^m}_{\cong} \\
   \hat{L}_\infty^\flat\otimes_{\varphi_q^m, {\mathbf{A}}_L}  \mathbf{A}^{H} \ar[r]^-{\id\varphi_q^m} & (\mathbb{C}_p^\flat)^H   }\]
shows that also the lower map is an isomorphism. Using that Verschiebung $V$ on  $ W_n((\mathbb{C}_p^\flat)^H)_L $ and $W_n(\hat{L}_\infty^\flat)_L$ is additive and satisfies the projection formula $V^m(x)\cdot y=V^m(x\cdot   \varphi_q^m(y))$ we see that we obtain a commutative exact diagram
\[\xymatrix{
  0  \ar[r]^{ } &  \hat{L}_\infty^\flat\otimes_{\varphi_q^n, {\mathbf{A}}_L}  \mathbf{A}^{H}  \ar[d]_{\id\varphi_q^n} \ar[r]^-{V^n \otimes \id} &  W_{n+1}(\hat{L}_\infty^\flat)_L\otimes_{ {\mathbf{A}}_L}  \mathbf{A}^{H} \ar[d]_{can} \ar[r]^{ } & W_n(\hat{L}_\infty^\flat)_L\otimes_{ {\mathbf{A}}_L}  \mathbf{A}^{H} \ar[d]_{\cong} \ar[r]^{ } & 0 \\
  0 \ar[r]^{ } & (\mathbb{C}_p^\flat)^H \ar[r]^-{V^n } & W_{n+1}((\mathbb{C}_p^\flat)^H)_L \ar[r]^{ } & W_n((\mathbb{C}_p^\flat)^H)_L ,     }\]
from which the claim follows by induction because the outer vertical maps are isomorphisms by the above and the induction hypothesis. Here the first non-trivial horizontal morphisms map onto the highest Witt vector component.

The second isomorphism is established as follows:  We choose a subgroup $N \subseteq H \subseteq H_L$ which is open normal in $H_L$ and obtain the extensions of henselian discrete valuation rings
\begin{equation*}
  \tilde{\mathbf{A}}^\dagger_L \subseteq (\tilde{\mathbf{A}}^\dagger)^H = W^\dagger((\mathbb{C}_p^\flat)^H)_L \subseteq (\tilde{\mathbf{A}}^\dagger)^N = W^\dagger((\mathbb{C}_p^\flat)^N)_L .
\end{equation*}
The corresponding extensions of their field of fractions
\begin{equation*}
  \tilde{\mathbf{B}}^\dagger_L \subseteq E := (\tilde{\mathbf{A}}^\dagger)^H [\tfrac{1}{\pi_L}]  \subseteq F := (\tilde{\mathbf{A}}^\dagger)^N [\tfrac{1}{\pi_L}]
\end{equation*}
satisfy $F^{H/N} = E$ and $F^{H_L/N} = \tilde{\mathbf{B}}^\dagger_L$. Hence $F/E$ and $F/\tilde{\mathbf{B}}^\dagger_L$ are Galois extensions of degree $[H:N]$ and $[H_L:N]$, respectively. It follows that $E/\tilde{\mathbf{B}}^\dagger_L$ is a finite extension of degree $[H_L:H]$. The henselian condition then implies\footnote{See Neukirch, Algebraische Zahlentheorie, proof of Satz II.6.8} that $(\tilde{\mathbf{A}}^\dagger)^H = W^\dagger((\mathbb{C}_p^\flat)^H)_L$ is free of rank $[H_L:H]$ over $\tilde{\mathbf{A}}^\dagger_L = W^\dagger(\hat{L}_\infty^\flat)_L$. The $\pi_L$-adic completion $(-)^{\widehat{}}$ of the two rings therefore can be obtained by the tensor product with $\tilde{\mathbf{A}}_L = W(\hat{L}_\infty^\flat)_L$. This gives the wanted
\begin{equation*}
  W(\hat{L}_\infty^\flat)_L \otimes_{\tilde{\mathbf{A}}^\dagger_L} (\tilde{\mathbf{A}}^\dagger)^H = W^\dagger(\hat{L}_\infty^\flat)_L^{\widehat{}} \otimes_{\tilde{\mathbf{A}}^\dagger_L} (\tilde{\mathbf{A}}^\dagger)^H = W^\dagger((\mathbb{C}_p^\flat)^H)_L^{\widehat{}} = W((\mathbb{C}_p^\flat)^H)_L .
\end{equation*}
\end{proof}

\begin{proposition}\label{f:exactint}
The   sequences
\begin{align}\label{f:exactA}
& 0 \rightarrow o_L \rightarrow {\mathbf{A}} \xrightarrow{\varphi_q - 1}{\mathbf{A}} \rightarrow 0,  \\\label{f:exacttildeA}
    & 0 \rightarrow o_L \rightarrow \tilde{\mathbf{A}} \xrightarrow{\varphi_q - 1}\tilde{\mathbf{A}} \rightarrow 0,  \\
    &  0 \rightarrow o_L \rightarrow\tilde{\mathbf{A}}^\dagger \xrightarrow{\varphi_q - 1} \tilde{\mathbf{A}}^\dagger \rightarrow 0.
\end{align}
are exact.
\end{proposition}

\begin{proof} The   first sequence is \cite[(26), Rem.\ 5.1]{SV15}. For the second sequence one proves by induction the statement for finite length Witt vectors using that the Artin-Schreier equation has a solution in $\mathbb{C}_p^\flat$. Taking projective limits then gives the claim. For the third sequence only the surjectivity has to be shown. This can be achieved by the same calculation as in the proof of \cite[Lem.\ 4.5.3]{KLII} with $R=\mathbb{C}_p^\flat$. \footnote{
 For the other see \cite[Lem.\ 4.5.3]{KLII} : There the exactness of  corresponding sequences for sheaves on the  pro\'{e}tale site $Spa(L,o_L)_{pro\acute{e}t}$ is shown, which in turn implies exactness for the corresponding sequences of stalks at the geometric point $Spa(\mathbb{C}_p, o_{\mathbb{C}_p})$. Note that taking stalks at this point is the same as taking sections over it.}
\end{proof}

\begin{lemma}\label{lem:ctssection}
For any finite $T$ in  $Rep_{o_L}(G_L)$ the map $\tilde{\mathbf{A}}\otimes_{o_L} T \xrightarrow{\varphi_q\otimes \id - 1} \tilde{\mathbf{A}}\otimes_{o_L} T$ has a continuous set theoretical section.
\end{lemma}

\begin{proof}
Since $T\cong \bigoplus_{i=1}^r o_L/\pi_L^{n_i}o_L$ for some natural numbers $r, n_i$ we may assume that $T=o_L/\pi_L^no_L$ for some $n$ and then we have to show that the surjective map $W_n(\mathbb{C}_p^\flat)_L \xrightarrow{\varphi_q-\id} W_n(\mathbb{C}_p^\flat)_L$ has a continuous set theoretical section. Thus me may neglect the additive structure and identify source and target with $X=(\mathbb{C}_p^\flat)^n$. In order to determine the components of the map $\varphi_q-\id =: f = (f_0,\ldots, f_{n-1}):X\to X$ with respect to these coordinates we recall that the addition in Witt rings is given by polynomials
\begin{equation*}\label{f:WittAdd}
  S_j(X_0,\ldots X_{j},Y_0,\ldots,Y_j)=X_j+Y_j + \mbox{ terms in }X_0,\ldots ,X_{j-1},Y_0,\ldots,Y_{j-1}
\end{equation*}
while the additive inverse is given by
\begin{equation*}\label{f:WittIn}
  I_j(X_0,\ldots X_{j} )= -X_j+  \mbox{ terms in }X_0,\ldots ,X_{j-1}.
\end{equation*}
Indeed, the polynomials $I_j$ are defined by the property that $\Phi_j(I_0,\ldots, I_j)=-\Phi_j(X_0,\ldots, X_j)$ where the Witt polynomials have the form $\Phi_j(X_0,\ldots,X_j)=X_0^{q^j}+\pi_LX_1^{q^{j-1}}+\ldots +\pi_L^jX_j.$ Modulo $(X_0,\ldots, X_{j-1})$ we derive that $\pi_L^j I_j(X_0,\ldots, X_j) \equiv -\pi_L^j X_j$ and the claim follows.
Since $\varphi_q$ acts componentwise rising the entries to their $q$th power, we conclude that
\[f_j=S_j(X_0^q,\ldots X_{j}^q,I_0(X_0),\ldots,I_j(X_0,\ldots X_{j} )).\] Hence the Jacobi matrix of $f$ at a point $x\in X$ looks like
\[D_x(f)=\begin{pmatrix}
           -1 &   & 0 \\
             &  \ddots &   \\
            \ast&   & -1 \\
         \end{pmatrix}
,\] i.e., is invertible in every point. As a polynomial map $f$ is locally analytic. It therefore follows from the inverse function theorem \cite[Prop.\ 6.4]{pLG} that $f$ restricts to a homeomorphism $f|U_0:U_0\xrightarrow{\cong}U_1$ of open neighbourhoods of $x$ and $f(x),$ respectively. By the surjectivity of $f$   every $x\in X$ has an open neighbourhood $U_x$ and a continuous map $s_x:U_x\to X$ with $f\circ s_x=\id_{|U_x}$. But $X$ is strictly paracompact by Remark 8.6 (i) in (loc.\ cit.), i.e., the covering $(U_x)_x$ has a disjoint refinement. There the restrictions of the $s_x$ glue to a continuous section of $f.$
\end{proof}

\begin{corollary}\label{cor:phi}
 For $T$ in  $Rep_{o_L}(G_L)$, the $n$th cohomology groups of the complexes concentrated  in degrees $0$ and $1$
     \begin{align}     \label{f:phitilde}
        &  \xymatrix@C=0.5cm{
          0 \ar[r] & \tilde{D}(T) \ar[rr]^{\varphi-1} && \tilde{D}(T) \ar[r] & 0 }\mbox{ and }\\
     \label{f:phi}
        &  \xymatrix@C=0.5cm{
          0 \ar[r] & {D}(T) \ar[rr]^{\varphi-1} && {D}(T) \ar[r] & 0   }
     \end{align}
     are isomorphic to $H^n(H_L,T)$ for any $n\geq 0$.
\end{corollary}

\begin{proof}
Assume first that $T$ is finite. For \eqref{f:phi} see \cite[Lemma 5.2]{SV15}. For \eqref{f:phitilde} we use Lemma \ref{lem:ctssection}, which says that the right hand map in the exact sequence
\[\xymatrix@C=0.5cm{
  0 \ar[r] & T \ar[rr]^{ } && \tilde{\mathbf{A}}\otimes_{o_L} T \ar[rr]^{\varphi_q\otimes \id -1} && \tilde{\mathbf{A}}\otimes_{o_L} T \ar[r] & 0 }\] has a continuous set theoretical section and thus gives rise to the long exact sequence of continuous cohomology groups
\begin{equation}\label{f:LES}
     0 \rightarrow H^0(H_L,T) \rightarrow \tilde{D}(T) \xrightarrow{\varphi-1} \tilde{D}(T) \rightarrow H^1(H_L,T) \rightarrow  H^1(H_L,\tilde{\mathbf{A}} \otimes_{o_L} T) \rightarrow \ldots
\end{equation}
Using the comparison isomorphism \eqref{f:comptildeD} and the subsequent Prop.\  \ref{prop:TatevanishingHII} we see that all terms from the fifth on vanish.

For the general case (for $\tilde{D}(T)$ as well as $D(T)$) we take inverse limits in the exact sequences for the $(T/\pi_L^mT)$ and observe that $H^n(H_L,T) \cong \varprojlim_{m} H^n(H_L,T/\pi_L^m T)$. This follows for $n \neq 2$ from \cite[Cor.\ 2.7.6]{NSW}. For $n=2$ we use \cite[Thm.\ 2.7.5]{NSW} and have to show that the projective system $(H^1(H_L,T/\pi^m_L T))_m$ is Mittag-Leffler. Since it is a quotient of the projective system $(D(T/\pi_L^m T))_m$, it suffices for this to check that the latter system is Mittag-Leffler. But due to the exactness of the functor $D$ this latter system is equal to the projective system of artinian $\mathbf{A}_L$-modules $(D(T)/\pi_L^m D(T))_m$ and hence is Mittag-Leffler. We conclude by observing that taking inverse limits of the system of sequences \eqref{f:LES} remains exact. The reasoning being the same for $\tilde{D}(T)$ and $D(T)$ we consider only the former. Indeed, we split the 4-term exact sequences into two short exact sequences of projective systems
\begin{equation*}
  0 \rightarrow H^0(H_L,V/\pi_L^m T) \rightarrow \tilde{D}(T/\pi_L^m T) \rightarrow (\varphi-1)\tilde{D}(T/\pi_L^m T) \rightarrow 0
\end{equation*}
and
\begin{equation*}
  0 \rightarrow (\varphi-1)\tilde{D}(T/\pi_L^m T) \rightarrow \tilde{D}(T/\pi_L^m T) \rightarrow H^1(H_L,T/\pi_L^m T) \rightarrow 0 .
\end{equation*}
Passing to the projective limits remains exact provided the left most projective systems have vanishing $\varprojlim^1$. For the system $H^0(H_L,T/\pi_L^m T)$ this is the case since it is Mittag-Leffler. The system $(\varphi-1)\tilde{D}(T/\pi_L^m T)$ even has surjective transition maps since the system $\tilde{D}(T/\pi_L^m T)$ has this property by the exactness of the functor $\tilde{D}$ (cf.\ Prop.\ \ref{Dbasechange}).
\end{proof}

\begin{proposition}\label{prop:TatevanishingHII}
$H^n(H,\tilde{\mathbf{A}}/\pi_L^m \tilde{\mathbf{A}})=0$ for all $n,m\geq 1$ and $H\subseteq H_L$  any closed subgroup.
 \end{proposition}

 \begin{proof}
 For $j<i$  the canonical projection $W_i(\mathbb{C}_p^\flat)\cong\tilde{\mathbf{A}}/\pi_L^i \tilde{\mathbf{A}}\twoheadrightarrow \tilde{\mathbf{A}}/\pi_L^j \tilde{\mathbf{A}}\cong W_j(\mathbb{C}_p^\flat)$ corresponds to the projection   $(\mathbb{C}_p^\flat)^i\twoheadrightarrow   (\mathbb{C}_p^\flat)^j $ and hence have set theoretical continuous sections. Using the associated long exact cohomology sequence (after adding the kernel) allows to reduce the statement to Prop.\  \ref{prop:TatevanishingH1}.
 \end{proof}

 For any commutative   ring $R$ with endomorphism $\varphi$ we write $\Phi(R)$ for the category of $\varphi$-modules consisting of $R$-modules equipped with a semi-linear $\varphi$-action.
We write $\Phi^{\acute{e}t}(R)$ for the subcategory of \'{e}tale $\varphi$-modules, i.e.,  such that $M$ is finitely generated over $R$ and $\varphi$ induces an $R$-linear isomorphism $\varphi^*M \xrightarrow{\cong} M$. Finally, we denote by   $\Phi^{\acute{e}t}_{f}(R)$ the subcategory consisting of finitely generated free $R$-modules.

For $M_1,M_2\in \Phi(R)$ the $R$-module $\Hom_R(M_1,M_2)$ has a natural structure as a $\varphi$-module satisfying
\begin{equation}\label{f:inner}
\varphi_{\Hom_{R}(M_1,M_2)}(\alpha)(\varphi_{M_1}(m))=\varphi_{M_2}( \alpha(m)) \ ,
\end{equation}
hence in particular
\begin{equation}\label{f:Homphiinv}
  \Hom_{R}(M_1,M_2)^{\varphi=\id}=\Hom_{\Phi(R)}(M_1,M_2).
\end{equation}
Note that with $M_1,M_2$ also $ \Hom_{R}(M_1,M_2)$ is \'{e}tale.
%If $M_1$ belongs to $\Phi^f(R)$, one has $\Hom_R(M_1,M_2)\cong M_1^\vee\otimes_R M_2$, where $M_1^\vee:=\Hom_R(M_1,R)$ denotes the dual $\varphi$-module.

\begin{remark}\label{rem:Yoneda}
 We recall from \cite[\S 1.5]{KLI} that the cohomology groups $H^i_\varphi(M)$ of the complex $M\xrightarrow{\varphi-1}M$ can be identified with the Yoneda extension groups $\mathrm{Ext}^i_{\Phi(R)}(R,M)$. Indeed, if $S:=R[X;\varphi]$ denotes the twisted polynomial ring satisfying $Xr=\varphi(r)X$ for all $r\in R,$ then we can identify $\Phi(R)$ with the category $S$-Mod of   (left) $S$-modules  by letting $X$ act via $\varphi_M$ on $X$. Using the free resolution
 \[\xymatrix@C=0.5cm{
   0 \ar[r] & S \ar[rr]^{\cdot(X-1)} && S \ar[rr]^{ } && R \ar[r] & 0 }\] the result follows.
\end{remark}

\begin{remark}\label{rem:faithfullyflat}
 Note that $\tilde{\bf A}_L^{\dagger}\subseteq \tilde{\bf A}_L$ is a faithfully flat ring extension as both rings are discrete valuation rings and the bigger one is the completion of the previous one.
\end{remark}

\begin{proposition}\label{prop:basechange} Base extension induces
\begin{enumerate}
  \item  an equivalence of categories
\[\Phi^{\acute{e}t}_{f}(\tilde{\bf A}_L^{\dagger} ) \leftrightarrow \Phi^{\acute{e}t}_{f}(\tilde{\bf A}_L^{})\]
  \item  and an isomorphism of Yoneda extension groups
  \[\mathrm{Ext}^1_{\Phi(\tilde{\bf A}_L^{\dagger} )}(\tilde{\bf A}_L^{\dagger},M) \cong \mathrm{Ext}^1_{\Phi(\tilde{\bf A}_L )}(\tilde{\bf A}_L, \tilde{\bf A}_L \otimes_{\tilde{\bf A}_L^{\dagger}} M) \] for all $M\in \Phi^{\acute{e}t}_{f}(\tilde{\bf A}_L^{\dagger} ).$
\end{enumerate}
\end{proposition}

\begin{proof}
  For the first item we imitate the proof of \cite[Thm.\ 8.5.3]{KLI}, see also \cite[Lem.\ 2.4.2,Thm.\ 2.4.5]{KedNew}: First we will show that for every $M\in \Phi^{\acute{e}t}_{f}(\tilde{\mathbf{A}}_L^\dagger)$ it holds that $(\tilde{\mathbf{A}}_L\otimes_{}M)^{\varphi=\id}\subseteq M^{\varphi=\id}$ and hence equality. Applied to $M:=\Hom_{\tilde{\mathbf{A}}_L^\dagger}(M_1,M_2)$ this implies that the  base change is  fully faithful by the equation \eqref{f:Homphiinv}. We observe that the analogue of \cite[Lem.\ 3.2.6]{KLI} holds in our setting and that $S$ in loc.\ cit.\ can be chosen to be a finite separable field extension of the perfect field $R=\hat{L}_\infty^\flat$. Thus we may choose $S$ in the analogue of \cite[Prop.\ 7.3.6]{KLI} (with $a=1$, $c=0$ and $M_0$ being our $M$) as completion of a (possibly infinite) separable field extension of $R$. This means in our situation that there exists a closed subgroup $H\subseteq H_L$ such that $(\tilde{\mathbf{A}}^\dagger)^H \otimes_{\tilde{\mathbf{A}}_L^\dagger} M
  = \bigoplus (\tilde{\mathbf{A}}^\dagger)^H e_i$ for a basis $e_i$ invariant under $\varphi.$ Now let $v=\sum x_i e_i$ be an arbitrary element in
  \begin{equation*}
    \tilde{\mathbf{A}}_L \otimes_{\tilde{\mathbf{A}}_L^\dagger}  M    \subseteq \tilde{\mathbf{A}}^H \otimes_{\tilde{\mathbf{A}}_L^\dagger} M
       = \tilde{\mathbf{A}}^H \otimes_{(\tilde{\mathbf{A}}^\dagger)^H}  (\tilde{\mathbf{A}}^\dagger)^H  \otimes_{\tilde{\mathbf{A}}_L^\dagger} M
       =\bigoplus \tilde{\mathbf{A}}^H e_i
  \end{equation*}
 with $x_i  \in \tilde{\mathbf{A}}^H$ and such that $\varphi(v)=v.$ The latter condition implies that $x_i\in \tilde{\mathbf{A}}^{H,\varphi_q=\id}=o_L,$ i.e., $v$ belongs to $(M\otimes_{\tilde{\mathbf{A}}_L^\dagger}(\tilde{\mathbf{A}}^\dagger)^H)\cap ( M\otimes_{\tilde{\mathbf{A}}_L^\dagger}\tilde{\mathbf{A}}_L)=M$, because $M$ is free and one has $\tilde{\mathbf{A}}_L\cap (\tilde{\mathbf{A}}^\dagger)^H=(\tilde{\mathbf{A}}^\dagger)^{H_L}= \tilde{\mathbf{A}}^\dagger_L.$
To show essential surjectivity one proceeds literally as in the proof of \cite[Thm.\ 8.5.3]{KLI} adapted to ramified Witt vectors. \Footnote{We may wish to discuss the equivalence among (a) and (d) in (loc.\ cit.), this is related to Prop \ref{prop:basechangephiGamma}.}

For the second statement choose a quasi-inverse functor $F: \Phi^{\acute{e}t}_{f}(\tilde{\bf A}_L^{})\to \Phi^{\acute{e}t}_{f}(\tilde{\bf A}_L^{\dagger} )$ with $F(\tilde{\bf A}_L)=\tilde{\bf A}_L^{\dagger}.$ Given an extension $\xymatrix@C=0.5cm{
    0 \ar[r] &  M\ar[rr]^{} && E \ar[rr]^{} && \tilde{\bf A}_L\ar[r] & 0 }$ over $\Phi(\tilde{\bf A}_L )$ with  $M\in \Phi^{\acute{e}t}_{f}(\tilde{\bf A}_L  )$ first observe that $E\in \Phi^{\acute{e}t}_{f}(\tilde{\bf A}_L )$, too. Indeed, $\tilde{\mathbf{A}}_L\xrightarrow{\varphi_q}\tilde{\mathbf{A}}_L$ is a flat ring extension, whence $\varphi^*E\to E$ is an isomorphism, if the corresponding outer maps are. The analogous statement holds over $\tilde{\bf A}_L^{\dagger}$. Therefore the sequence $\xymatrix@C=0.5cm{
    0 \ar[r] &  F(M)\ar[rr]^{} && F(E) \ar[rr]^{} && \tilde{\bf A}_L^{\dagger}\ar[r] & 0 }$ is exact by Remark \ref{rem:faithfullyflat},  because its base extension - being isomorphic to the original extension - is, by assumption.
\end{proof}

%\begin{proposition}
% The canonical map
% \begin{equation}\label{f:Ddaggerbasechange}
%    \tilde{\mathbf{A}}_L\otimes_{\tilde{ {\mathbf{A}}}_L^\dagger} \tilde{{D}}^\dagger(T)\xrightarrow{\cong} \tilde{D}(V)
% \end{equation} is an isomorphism.
% %and
%% the functor   $\tilde{D}(-):Rep_{o_L}(G_L) \to \mathfrak{M}(\tilde{\mathbf{A}}_L)$ is exact. Moreover, we have a comparison isomorphisms
%%     \begin{align}\label{f:comptildeD}
%%      \tilde{\mathbf{A}}\otimes_{o_L} T & \cong  \tilde{\mathbf{A}}\otimes_{ \tilde{\mathbf{A}}_L} \tilde{D}(T).
%%     \end{align}
%\end{proposition}
%
%\begin{proof}
%The proof follows the same strategy as that of \eqref{f:Dbasechange}.
%\end{proof}

We denote by $\mathfrak{M}^{\acute{e}t}_{f}(\tilde{\bf A}_L^{\dagger} )$ and $\mathfrak{M}^{\acute{e}t}_{f}(\tilde{\bf A}_L^{})$ the full subcategories of $\mathfrak{M}^{\acute{e}t}(\tilde{\bf A}_L^{\dagger} )$ and $\mathfrak{M}^{\acute{e}t}(\tilde{\bf A}_L^{})$, respectively, consisting of finitely generated free modules over the base ring.

\begin{remark}\label{rem:inducedCanonicalTop}
Let $M$ be in $  \mathfrak{M}^{\acute{e}t}_{f}(\tilde{\bf A}_L^{})$ and endow $N:=\tilde{\bf A}_L\otimes_{\tilde{\bf A}_L^{\dagger}} M$ with the canonical topology with respect to the weak topology of $\tilde{\bf A}_L $. Then the induced subspace topology of $M\subseteq N$ coincides with the canonical topology with respect to the weak topology of $\tilde{\bf A}_L^{\dagger} $. Indeed for free modules this is obvious while for torsion modules this can be reduced by the elementary divisor theory to the case $M=\tilde{\bf A}_L^{\dagger}/\pi_L^n \tilde{\bf A}_L^{\dagger}\cong \tilde{\bf A}_L/\pi_L^n \tilde{\bf A}_L$. But the latter spaces are direct product factors of $\tilde{\bf A}_L^{\dagger}$ and $\tilde{\bf A}_L,$ respectively, as topological spaces, from wich the claim easily follows.
\end{remark}

\begin{proposition}\label{prop:basechangephiGamma}
For $T\in\Rep_{o_L}(G_L)$ and $V\in\Rep_L(G_L)$  we have  natural isomorphisms
\begin{align}\label{f:basechangedagger2}
\tilde{\bf A}_L^{}\otimes_{\tilde{\bf A}_L^{\dagger}}\tilde{D}^\dagger(T) & \cong \tilde{D}(T)\mbox{ and }\\
  \label{f:basechangedagger3}\tilde{\bf B}_L^{}\otimes_{\tilde{\bf B}_L^{\dagger}}\tilde{D}^\dagger(V) & \cong \tilde{D}(V),
\end{align}
as well as
\begin{align}\label{f:compdagger2}
\tilde{\bf A}^{\dagger}\otimes_{\tilde{\bf A}_L^{\dagger}}\tilde{D}^\dagger(T) & \cong \tilde{\bf A}^{\dagger}\otimes_{o_L} T\mbox{ and }\\
  \label{f:compdagger3}\tilde{\bf B}^{\dagger}\otimes_{\tilde{\bf B}_L^{\dagger}}\tilde{D}^\dagger(V) & \cong \tilde{\bf B}^{\dagger}\otimes_L V,
\end{align}
respectively.  In particular,  the functor  $\tilde{D}^\dagger(-):\Rep_{o_L}(G_L) \to \mathfrak{M}^{\acute{e}t}(\tilde{\mathbf{A}}_L^\dagger)$ is exact.

Moreover, base extension induces  equivalences of categories
\[\mathfrak{M}^{\acute{e}t}_f(\tilde{\bf A}_L^{\dagger} ) \leftrightarrow \mathfrak{M}^{\acute{e}t}_f(\tilde{\bf A}_L^{}),\]
and hence also an equivalence of categories
\[\mathfrak{M}^{\acute{e}t}(\tilde{\bf B}_L^{\dagger} ) \leftrightarrow \mathfrak{M}^{\acute{e}t}(\tilde{\bf B}_L^{}).\]
\end{proposition}

\begin{proof} Note that  the base change functor is well-defined - regarding the continuity of the $\Gamma_L$-action - by Lemma \ref{lem:contGammaAction} and Remark \ref{rem:faithfullyflat}  while $\tilde{D}^\dagger$ is well-defined by Remark \ref{rem:inducedCanonicalTop}, once \eqref{f:basechangedagger2} will have been shown. We first show the equivalence of categories for free modules:
By   Prop.\  \ref{prop:basechange} we already have, for $M_1,M_2\in\mathfrak{M}^{\acute{e}t}_{f}(\tilde{\bf A}_L^{\dagger} )$,  an isomorphism
\[
\Hom_{\Phi(\tilde{\bf A}_L^{\dagger} )} (M_1,M_2) \cong \Hom_{\Phi(\tilde{\bf A}_L )} (\tilde{\bf A}_L \otimes_{\tilde{\bf A}_L^{\dagger}} M_1,\tilde{\bf A}_L \otimes_{\tilde{\bf A}_L^{\dagger}} M_2).
\]
Taking $\Gamma_L$-invariants gives that the base change functor in question is fully faithful.

In order to show that this base change functor is also essentially surjective,
\Footnote{Alternatively imitate proof of \cite[Thm.\ 2.4.5]{KedNew} perhaps, but be careful as the proof of \cite[Thm.\ 2.6.2]{KedNew} cannot be true as otherwise all reps would be overconvergent (also in the imperfect case) - vielleicht geht ja Lemma 2.5.1 in (loc.cit.) schief!?} consider an arbitrary $N \in \mathfrak{M}^{\acute{e}t}_{f}(\tilde{\bf A}_L^{}).$ Again by  \ref{prop:basechange} we know that there is  a free \'{e}tale $\varphi$-module $M$ over $\tilde{\bf A}_L^{\dagger} $ whose base change is isomorphic to $N$. By the fully faithfulness the $\Gamma_L$-action  descends to  $M$\footnote{As $\gamma\in\Gamma_L$ acts semilinearly, one formally has to replace $N\xrightarrow{\gamma} N$ by the linearized isomorphism $\tilde{\bf A}_L\otimes_{\gamma,\tilde{\bf A}_L}N\xrightarrow{\gamma^{lin}} N$. Upon checking that the source is again a \'{e}tale $\varphi$-module with model $\tilde{\bf A}_L^\dagger\otimes_{\gamma,\tilde{\bf A}_L^\dagger}M$ one sees by the fully faithfulness on $\varphi$-modules that the linearized isomorphism descends and induces the desired semi-linear action.}.   Since the weak topology of $M$ is compatible with that of $N$ by Remark \ref{rem:inducedCanonicalTop}, this action is again continuous.

To prepare for the proof of the isomorphism \eqref{f:basechangedagger2} we first observe the following fact. The isomorphism \eqref{f:comptildeD} implies that $T$ and $\tilde{D}(T)$ have the same elementary divisors, i.e.: If $T \cong \oplus_{i=1}^r o_L/\pi_L^{n_i} o_L$ as $o_L$-module (with $n_i \in \mathbb{N} \cup \{\infty\}$) then $\tilde{D}(T) \cong \oplus_{i=1}^r \tilde{\mathbf{A}}_L / \pi_L^{n_i} \tilde{\mathbf{A}}_L$ as $\tilde{\mathbf{A}}_L$-module.

We shall prove \eqref{f:basechangedagger2} in several steps: First assume that $T$ is \textit{finite}. Then $T$ is annihilated by some $\pi_L^n$. We have $\tilde{D}^\dagger(T)= \tilde{D}(T)$ and $\tilde{\mathbf{A}}_L^\dagger / \pi_L^n \tilde{\mathbf{A}}_L^\dagger = \tilde{\mathbf{A}}_L / \pi_L^n \tilde{\mathbf{A}}_L$ so that there is nothing to prove. Secondly we suppose that $T$ is \textit{free} {\it and}  that $\tilde{D}^\dagger(T)$ is free over $\tilde{\mathbf{A}}_L^\dagger$ of the same rank $r := \rk_{o_L} T$. On the other hand, as the functor $\tilde{D}^\dagger$ is always left exact, we obtain the injective maps
\[
     \tilde{D}^\dagger(T) / \pi_L^n  \tilde{D}^\dagger(T) \rightarrow \tilde{D}^\dagger(T/\pi_L^n T) = \tilde{D}(T/ \pi_L^n T) .
\]
for any $n \geq 1$. We observe that both sides are isomorphic to $(\tilde{\mathbf{A}}_L^\dagger / \pi_L^n \tilde{\mathbf{A}}_L^\dagger)^r = (\tilde{\mathbf{A}}_L / \pi_L^n \tilde{\mathbf{A}}_L)^r$. Hence the above injective maps are bijections. We deduce that
\begin{align*}
\tilde{\mathbf{A}}_L\otimes_{ {\mathbf{A}}_L^\dagger} \tilde{D}^\dagger(T)   & \cong  \varprojlim_{n} \tilde{D}^\dagger(T)/\pi_L^n\tilde{D}^\dagger(T)  \\
    &  \cong   \varprojlim_{n} \tilde{D}(T /\pi_L^n T)  \\
    &  \cong    \varprojlim_{n} \tilde{D}(T)/\pi_L^n\tilde{D}(T) \\
    &  \cong  \tilde{D}(T)
\end{align*}
using that the above tensor product means $\pi_L$-adic completion  for finitely generated $\tilde{\mathbf{A}}^\dagger_L$-modules.

Thirdly  let $T \in \Rep_{o_L,f}(G_L)$ be arbitrary and  $M \in \mathfrak{M}^{\acute{e}t}_{f}(\tilde{\bf A}_L^{\dagger} )$ such that $\tilde{\bf A}_L \otimes_{\tilde{\bf A}_L^\dagger} M \cong \tilde{D}(T)$ according the equivalence of categories.    Without loss of generality we may treat this isomorphism as an equality. Similarly as in the proof of Prop.\  \ref{prop:basechange} and with the same notation one shows that $(\tilde{\mathbf{A}}^\dagger \otimes_{\tilde{\mathbf{A}}^\dagger_L } M)^{\varphi =1} = \bigoplus_{i=1}^ro_Le_i$ for some appropriate $\varphi$-invariant basis $e_1,\ldots, e_r$ of $\tilde{\mathbf{A}}^\dagger \otimes_{\tilde{\mathbf{A}}^\dagger_L } M$. Note that $r = \rk_{o_L} T$. Using \eqref{f:comptildeD}, it follows that
\begin{align*}
  T & = (\tilde{\mathbf{A}} \otimes_{o_L} T)^{\varphi=1} \cong ( \tilde{\mathbf{A}} \otimes_{\tilde{\mathbf{A}}_L} \tilde{D}(T))^{\varphi=1} = (\tilde{\mathbf{A}} \otimes_{\tilde{\mathbf{A}}_L^\dagger } M)^{\varphi=1}    \\
  & = \bigoplus_{i=1}^r\tilde{\mathbf{A}}^{\varphi_q=1} e_i = \bigoplus_{i=1}^r o_Le_i = (\tilde{\mathbf{A}}^\dagger \otimes_{\tilde{\mathbf{A}}^\dagger_L } M)^{\varphi=1} .
\end{align*}
It shows that the comparison isomorphism \eqref{f:comptildeD} restricts to an injective map $T \hookrightarrow \tilde{\mathbf{A}}^\dagger \otimes_{\tilde{\mathbf{A}}^\dagger_L } M$, which extends to a homomorphism $\tilde{\mathbf{A}}^\dagger \otimes_{o_L} T \xrightarrow{\alpha} \tilde{\mathbf{A}}^\dagger \otimes_{\tilde{\mathbf{A}}^\dagger_L } M$ of free $\tilde{\mathbf{A}}^\dagger$-modules of the same rank $r$. Further base extension by $\tilde{\mathbf{A}}$ gives back the isomorphism \eqref{f:comptildeD}. Since $\tilde{\mathbf{A}}$ is faithfully flat over $\tilde{\mathbf{A}}^\dagger$ the map $\alpha$ was an isomorphism already. By passing to $H_L$-invariants we obtain an isomorphism $\tilde{D}^\dagger(T) \cong M$ and see that $\tilde{D}^\dagger(T)$ is free of the same rank as $T$. Hence the second case applies and gives \eqref{f:basechangedagger2} for free $T$ and \eqref{f:basechangedagger3}. Finally, let $T$ be just finitely generated over $o_L.$ Write $0\to T_{\mathrm{fin}}\to T\to T_{\mathrm{free}}\to 0$ with finite $T_{\mathrm{fin}}$  and free  $T_{\mathrm{free}}.$ We then have the commutative exact diagram
\[\small\xymatrix@C=0.5cm{
  0  \ar[r]^{ } & \tilde{\mathbf{A}}_L\otimes_{\tilde{\mathbf{A}}_L^\dagger} \tilde{D}^\dagger(T_{\mathrm{fin}}) \ar[d]_{ \cong} \ar[r]^{ } &  \tilde{\mathbf{A}}_L\otimes_{\tilde{\mathbf{A}}_L^\dagger} \tilde{D}^\dagger(T) \ar[d]_{} \ar[r]^{} & \tilde{\mathbf{A}}_L\otimes_{\tilde{\mathbf{A}_L^\dagger}} \tilde{D}^\dagger(T_{\mathrm{free}}) \ar[d]_{\cong} \ar[r]^{ } &  \tilde{\mathbf{A}}_L\otimes_{\tilde{\mathbf{A}}_L^\dagger} H^1(H_L,\tilde{\mathbf{A}}^\dagger\otimes_{o_L} T_{\mathrm{fin}}) \\
  0 \ar[r]^{ } &   \tilde{D}(T_{\mathrm{fin}})\ar[r]^{ } & \tilde{D}(T)\ar[r]^{ } &  \tilde{D}(T_{\mathrm{free}}) \ar[r]^{} & 0,   }\]
  in which we use the first and third step for the vertical isomorphisms.
In order to show that the middle perpendicular arrow is an isomorphism it suffices to prove that $  H^1(H_L,\tilde{\mathbf{A}}^\dagger\otimes_{o_L} T_{\mathrm{fin}}) =0$. But since $T_{\mathrm{fin}}$ is annihilated by some $\pi_L^n$ we have
\[\tilde{\mathbf{A}}^\dagger\otimes_{o_L} T_{\mathrm{fin}}\cong \tilde{\mathbf{A}}/\pi_L^n\tilde{\mathbf{A}}\otimes_{o_L} T_{\mathrm{fin}}\cong\tilde{\mathbf{A}}/\pi_L^n\tilde{\mathbf{A}}\otimes_{\tilde{\mathbf{A}}_L}\tilde{D}( T_{\mathrm{fin}}),\]
the last isomorphism by \eqref{f:comptildeD}. Thus it suffices to prove the vanishing of $  H^1(H_L,\tilde{\mathbf{A}}/\pi_L^n\tilde{\mathbf{A}}),$ which is established in Prop.\  \ref{prop:TatevanishingHII} and finishes the  proof of the isomorphism \eqref{f:basechangedagger2}.

Note that this base change isomorphism implies the exactness of $\tilde{D}^\dagger$ as $\tilde{D} $ is exact by Prop.\  \ref{Dbasechange} and using that the base extension is faithfully flat by Remark \ref{rem:faithfullyflat}.

For free $T$ the statement \eqref{f:compdagger2} (and hence \eqref{f:compdagger3}) is already implicit in the above arguments while for finite $T$ the statement coincides with \eqref{f:comptildeD}. The general case follows from the previous ones by exactness of $\tilde{D}^\dagger$ and the five lemma as above.

\end{proof}

\begin{corollary}\label{cor:phitildedagger}
 For a   $T$ in  $\Rep_{o_L,f}(G_L)$ and $V$ in  $\Rep_{L}(G_L)$, the $n$th cohomology group, for any $n \geq 0$, of the complexes concentrated  in degrees $0$ and $1$
     \begin{align}     \label{f:phitildedagger}
        &  \xymatrix@C=0.5cm{
          0 \ar[r] & \tilde{D}^\dagger(T) \ar[rr]^{\varphi-1} && \tilde{D}^\dagger(T) \ar[r] & 0 }\mbox{ and } \\
           \label{f:phitildedaggerV}
        &  \xymatrix@C=0.5cm{
          0 \ar[r] & \tilde{D}^\dagger(V) \ar[rr]^{\varphi-1} && \tilde{D}^\dagger(V) \ar[r] & 0 }\mbox{ and }
     \end{align}
     is isomorphic to $H^n(H_L,T)$  and $H^n(H_L,V),$  respectively.
\end{corollary}

\begin{proof}
 The integral result reduces, by \eqref{f:basechangedagger2}, Remark \ref{rem:Yoneda}, and Prop.\  \ref{prop:basechange}, to Corollary \ref{cor:phi}. Since inverting $\pi_L$ is exact and commutes with taking cohomology \cite[Prop.\ 2.7.11]{NSW}, the second statement follows.
\end{proof}

Set ${\bf A}^\dagger:=\tilde{\bf A}^{\dagger}\cap {\bf A}$ and ${\bf B}^\dagger:={\bf A}^\dagger[\frac{1}{\pi_L}]$ as well as ${\bf A}^\dagger_L:=({\bf A}^\dagger)^{H_L}.$ Note that ${\bf B}^\dagger_L:=({\bf B}^\dagger)^{H_L}\subseteq {\bf B}^\dagger\subseteq{\tilde{\bf B}}^\dagger.$
For $V\in \Rep_L(G_L)$ we define $D^\dagger(V):=(\mathbf{B}^\dagger\otimes_L V)^{H_L}.$ The categories $\mathfrak{M}^{\acute{e}t}({\bf A}^\dagger_L)$ and $ \mathfrak{M}^{\acute{e}t}({\bf B}^\dagger_L)$ are defined analogously as in Definition \ref{def:etale}.

\emph{\begin{remark}\label{rem:overconvergent}
 There  is also the following more concrete description for ${\bf A}^\dagger_L$ in terms of Laurent series in $\omega_{LT}:$
\begin{align*}
{\bf A}^\dagger_L=\{F(\omega_{LT})\in \mathbf{A}_L|& F(Z)\mbox{ converges on } \rho\leq |Z|<1 \mbox{ for some } \rho\in (0,1)\}\subseteq \mathbf{A}_L.
\end{align*}
Indeed this follows from  the analogue of \cite[Lem.\ II.2.2]{ChCo} upon noting that the latter holds with and without the integrality condition:  ''$rv_p(a_n)+n\geq 0$ for all $n\in\mathbb{Z}$'' (for $r\in \overline{\mathbf{R}}\setminus \mathbf{R}$) in the notation of that article.
\Footnote{ This description does not require any completeness property! A similar result holds for $\mathbf{A}_{\leq 1,L}^\dagger$ when requiring for the Laurent series in $\omega_{LT}$ in addition that $F(Z)$ takes values on $ \rho\leq |Z|<1$ of norm at most $1.$ More precisely, for $r<1$ (or equivalently $s(r):=\frac{q-1}{rq}>\frac{q-1}{q}$) $W^r(\mathbf{C}_p^\flat)_L$ and $W^r_{\leq1}(\mathbf{C}_p^\flat)_L$ correspond to  $\{F(\omega_{LT})\in\mathbf{A}_L| \; F(Z) \mbox{ converges on } |\omega|^r\leq |Z|<1\}$ and $\{F(\omega_{LT})\in\mathbf{A}_L| \; F(Z) \mbox{ converges on } |\omega|^r=q^{-\frac{1}{s(r)}}\leq |Z|<1 \mbox{ with values } |F(z)|\leq1 \}$, respectively. The latter condition on the values can also be rephrased as $s(r)v_{\pi_L}(a_m)+m\geq 0$ for all $m\in \z$ corresponding to $V(x,s(r))\geq 0$ on the Witt vector side if $F(Z)=\sum_m a_mZ^m\in\mathbf{A}_L.$
}
In particular we obtain canonical embeddings ${\bf A}^\dagger_L\subseteq {\bf B}^\dagger_L \hookrightarrow \cR_L$ of rings.
\end{remark}}

\begin{definition}\label{def:overconv}
$V$ in $\Rep_L(G_L)$ is called {\it overconvergent}, if $\dim_{{\bf B}^\dagger_L} D^\dagger(V) = \dim_LV.$  We denote by   $\mathrm{Rep}_L^\dagger(G_{L})\subseteq \Rep_L(G_L)$ the full subcategory of overconvergent representations.
\end{definition}

\begin{remark}
   We always have $\dim_{{\bf B}^\dagger_L} D^\dagger(V) \leq \dim_LV$. If $V \in \Rep_L(G_L)$ is overconvergent then we have the natural isomorphism
\begin{equation}\label{f:overconv}
  {\mathbf{B}}_L \otimes_{{\mathbf{B}}^\dagger_L}D^\dagger(V) \xrightarrow{\cong} {D}(V).
\end{equation}
\end{remark}
\begin{proof}
Since $\mathbf{B}_L$ and $\mathbf{B}^\dagger_L$ are fields this is immediate from \cite[Thm.\ 2.13]{FO}.
\end{proof}

\begin{remark}
In \cite[\S 10]{Be16} Berger uses the following condition to define overconvergence of $V$: There exists a $\mathbf{B}_L$-basis $x_1,\ldots,x_n$ of $D(V)$ such that $M := \bigoplus_{i=1}^n \mathbf{B}^\dagger_L x_i$ is a $(\varphi_L,\Gamma_L)$-module over  $\mathbf{B}_L^\dagger.$ This then implies  a natural isomorphism
\begin{equation}\label{f:overconBerger}
  \mathbf{B}_L \otimes_{\mathbf{B}^\dagger_L } M \cong D(V).
\end{equation}
\end{remark}

\begin{lemma}
$V$ in $\Rep_L(G_L)$ is overconvergent if and only if $V$ satisfies the above condition of Berger. In this case $M=D^\dagger(V).$
\end{lemma}

\begin{proof}
If $V$ is overconvergent, we can take a basis within $M:=D^\dagger(V).$ Conversely let $V$ satisfy Berger's condition, i.e.\ we have the isomorphism \eqref{f:overconBerger}. One easily checks by faithfully flat descent that with $D(V)$ also $M$ is \'{e}tale. By \cite[Prop.\ 1.5 (a)]{FX}\footnote{Note that there $\bar{D}$ actually belongs to the category of $(\varphi,G_F)$-modules over $\tilde{\mathbf{B}}_{\mathbb{Q}_p}\otimes F $ instead of over $\tilde{\mathbf{B}}_{\mathbb{Q}_p} $ in their notation. } we obtain the identity $V=\left( \mathbf{B}^\dagger\otimes_{\mathbf{B}^\dagger_L} M\right)^{\varphi=1}$ induced from the comparison isomorphism
\begin{equation}\label{f:compBergerover}
\mathbf{B}\otimes_L V\cong \mathbf{B}\otimes_{\mathbf{B}_L}D(V)\cong \mathbf{B}\otimes_{\mathbf{B}_L^\dagger}M.
\end{equation}
 We shall prove that $M\subseteq D^\dagger(V)=(\mathbf{B} ^\dagger\otimes_LV)^{H_L}$ as then $M=D^\dagger(V)$ by dimension reasons. To this aim we may write a basis $v_1,\ldots, v_n$ of $V$ over $L$ as $v_i=\sum c_{ij}x_j$ with $c_{ij}\in \mathbf{B}^\dagger$. Then \eqref{f:compBergerover} implies that the matrix $C=(c_{ij})$ belongs to $M_n(\mathbf{B}^\dagger )\cap GL_n( \mathbf{B} )=GL_n(\mathbf{B}^\dagger).$ Thus $M$ is contained in $\mathbf{B}^\dagger\otimes_LV$ and - as subspace of $D(V)$ - also $H_L$-invariant, whence the claim.
\end{proof}

\begin{remark}\label{rem:notoverconvergent}
Note that the imperfect version of Prop.\   \ref{prop:basechangephiGamma} is not true: the base change $\mathfrak{M}^{\acute{e}t}({\mathbf{B}}_L^\dagger) \to \mathfrak{M}^{\acute{e}t}({\mathbf{B}}_L) $ is not essentially surjective in general, whence not an equivalence of categories, by \cite{FX}.  By definition, its essential image consists of {\it overconvergent} $(\varphi_L,\Gamma_L)$-modules, i.e., whose  corresponding Galois representations are   {\it overconvergent}.
\end{remark}

\begin{lemma} Assume that $V\in \Rep_L(G_L)$ is overconvergent. Then there is natural isomorphism
\[\tilde{\mathbf{B}}^\dagger_L\otimes_{\tilde{\mathbf{B}}^\dagger_L}D^\dagger(V)\cong \tilde{D}^\dagger(V).\]
\end{lemma}

\begin{proof}
By construction we have a natural map $\tilde{\mathbf{B}}^\dagger_L\otimes_{\tilde{\mathbf{B}}^\dagger_L}D^\dagger(V)\to\tilde{D}^\dagger(V),$ whose base change to $\tilde{\mathbf{B}}_L$
\[\tilde{\mathbf{B}}_L\otimes_{\tilde{\mathbf{B}}^\dagger_L}D^\dagger(V)\to\tilde{\mathbf{B}}_L\otimes_{\tilde{\mathbf{B}}^\dagger_L}\tilde{D}^\dagger(V)\cong \tilde{D}(V)  \]
arises also as the base change of the isomorphism \eqref{f:overconv}, %${\mathbf{B}}_L\otimes_{{\mathbf{B}}^\dagger_L}D^\dagger(V)\cong{D}(V)$,
 whence is an isomorphism itself. Here we have used the (base change of the) isomorphisms \eqref{f:basechangedagger3}, \eqref{f:Dbasechange}. By faithfully flatness the  original map is an isomorphism, too.
\end{proof}

\section{The perfect Robba ring}

Again let $K$ be   any perfectoid field containing $L$ and  $r>0.$ For $0<s\leq r,$ let $\tilde{\cR}^{[s,r]}(K)$ be the completion of  $W^r(K^\flat)_L[\frac{1}{\pi_L}]$ with respect to the norm $\max\{|\ |_{s} , |\ |_{r}\},$ and put
%\footnote{ The $W^r(K^\flat)_L[\frac{1}{\pi_L}]$ are normed $L$-vector spaces with respect to  $|\ |_{s}.$ The continuity of the transition maps in the inverse limit below is easily checked using the norms $|\ |_{s}^\frac{1}{s},$ which are NOT $L$-vector space norms!}
\[
\tilde{\cR}^r(K) = \varprojlim_{s\in (0,r]} \tilde{\cR}^{[s,r]}(K)
\]
equipped with the Fr\'{e}chet topology. Let $\tilde{\cR}(K) = \varinjlim_{r>0} \tilde{\cR}^r(K)$, equipped with the locally convex direct limit topology (LF topology). We set
$\tilde{\cR}=\tilde{\cR}({\mathbb{C}_p})$ and  $\tilde{\cR}_{L}:=\tilde{\cR}({\hat{L}_\infty})$. For  geometric interpretation  of these definitions, see \cite{Ede}. As in \cite[Thm.\ 9.2.15]{KLI} we have
\[\tilde{\cR}^{H_L}=\tilde{\cR}_{L}.\]
Recall from section \ref{sec:notation} the embedding  $o_L[[Z]]\to  W(\tilde{\mathbf{E}})_L$. As we will explain in section \ref{sec:weak} the image $\omega_{LT}$ of the variable $Z$ already lies in $W({\hat{L}_\infty}^\flat)_L$, so that we actually have an embedding $o_L[[Z]]\to  W({\hat{L}_\infty}^\flat)_L$.
Similarly as in \cite[Def.\ 4.3.1]{KLI} for the cyclotomic situation one shows that the latter  embedding   extends to a $\Gamma_L$- and $\varphi_L$-equivariant topological monomorphism  $\mathcal{R}_L\to \tilde{\mathcal{R}}_L,$ see also \cite[Konstruktion 1.3.27]{W} in the Lubin-Tate setting.

Let $R$ be either $\mathcal{R}_L$ or $\tilde{\mathcal{R}}_L$. A $(\varphi_L,\Gamma_L)$-module over $R$ is a finitely generated free $R$-module $M$ equipped with commuting semilinear actions of $\varphi_M$ and $\Gamma_L$, such that the action is continuous for  the LF topology and such that the semi-linear map $\varphi_M:M\to M$ induces an isomorphism $\varphi_M^{lin}:R\otimes_{R,\varphi_R}M\xrightarrow{\cong} M.$ Such $M$ is called \'{e}tale, if there exists an \'{e}tale $(\varphi_L,\Gamma_L)$-module $N$ over $\mathbf{A}_L^\dagger$ and $\tilde{\mathbf{A}}_L^\dagger$ (see before Definition \ref{def:etale}), such that $\mathcal{R}_L\otimes_{\mathbf{A}_L^\dagger} N\cong M$ and $\tilde{\mathcal{R}}_L\otimes_{\tilde{\mathbf{A}}_L^\dagger} N\cong M,$ respectively.
%\footnote{\label{foot:KedlayaSlope} This is compatible with our earlier
%Definition \ref{def:etaleana} \cite[Chap.\ 4]{KLI}.}

By $\mathfrak{M}(R)$ and  $\mathfrak{M}^{\acute{e}t}(R)$ we denote the category of $(\varphi_L,\Gamma_L)$-modules and \'{e}tale $(\varphi_L,\Gamma_L)$-modules over $R,$ respectively.

We call the topologies on $\tilde{\mathbf{A}}_L^\dagger$  and $\tilde{\mathbf{A}}^\dagger$, which make the inclusions $\tilde{\mathbf{A}}_L^\dagger\subseteq \tilde{\mathbf{A}}^\dagger\subseteq\tilde{\cR}$ topological embeddings,  the LF-topologies.

\begin{lemma}\label{lem:LF}
For $M\in \mathfrak{M}^{\acute{e}t}_f(\tilde{\mathbf{A}}_L^\dagger)$ the $\Gamma_L$-action is also continuous with respect to the canonical topology with respect to the LF-topology of $\tilde{\mathbf{A}}_L^\dagger$.
\end{lemma}

\begin{proof} The proof in fact works in the following generality: Suppose that $\tilde{\mathbf{A}}^\dagger$ is equipped with an $o_L$-linear ring topology which induces the $\pi_L$-adic topology on $o_L$. Consider on $\tilde{\mathbf{A}}^\dagger_L$ the corresponding induced topology. We claim that then the $\Gamma_L$-action on $M$ is continuous with respect to the corresponding canonical topology.
 By Prop.\  \ref{prop:equivalenceintegral} we may choose  $T\in\Rep_{o_L,f}(G_L)$ such that $M\cong \tilde{D}^\dagger(T).$ Then we have a homeomorphism $\tilde{\mathbf{A}}^\dagger\otimes_{o_L}T\cong\tilde{\mathbf{A}}^\dagger\otimes_{\tilde{\mathbf{A}}^\dagger_L}M$ with respect to the canonical topology by \eqref{f:compdagger2} (as any $R$-module homomorphism of finitely generated modules is continuous with respect to the canonical topology with regard to any topological ring $R$). Since $o_L\subseteq \tilde{\mathbf{A}}^\dagger$ is a topological embedding with respect to the $\pi_L$-adic and the given topology, respectively, Lemma \ref{lem:contGammaAction} implies that $G_L$ is acting continuously on $\tilde{\mathbf{A}}^\dagger\otimes_{\tilde{\mathbf{A}}^\dagger_L}M$, whence $\Gamma_L$ acts continuously on $M=\left( \tilde{\mathbf{A}}^\dagger\otimes_{\tilde{\mathbf{A}}^\dagger_L}M\right)^{H_L}$ with respect to the induced topology as subspace of the previous module. Since all involved modules are free and hence carry the product topologies and since $\tilde{\mathbf{A}}^\dagger_L\subseteq \tilde{\mathbf{A}}^\dagger$ is a topological embedding, it is clear that the latter topology of $M$ coincides with its canonical topology.
\end{proof}

We define the functor
\begin{align*}
  \tilde{D}^\dagger_{rig}(-) : \Rep_{L}(G_L) & \longrightarrow \mathfrak{M}(\tilde{\mathcal{R}}_L) \\
                           V & \longmapsto (\tilde{\mathcal{R}} \otimes_{L} V)^{H_L} ,
\end{align*}
 where the fact, that $\Gamma_L$ acts continuously on the image with respect to the LF-topology can be seen as follows, once we have shown the next lemma. Indeed, \eqref{f:Dbasechangedaggerrig} implies that for any $G_L$-stable $o_L$-lattice $T$ of $V$ we also have an isomorphism $ \tilde{\mathcal{R}}_L \otimes_{\tilde{\mathbf{A}}^\dagger_L} \tilde{D}^\dagger(T)\xrightarrow{\cong} \tilde{D}^\dagger_{rig}$. Now   again  Lemma \ref{lem:contGammaAction} applies to conclude the claim.

\begin{lemma}\label{lem:daggerrig}
 The canonical map
 \begin{equation}\label{f:Dbasechangedaggerrig}
  \tilde{\mathcal{R}}_L \otimes_{\tilde{\mathbf{B}}^\dagger_L} \tilde{D}^\dagger(V)\xrightarrow{\cong} \tilde{D}^\dagger_{rig}(V)
 \end{equation} is an isomorphism and
 the functor   $ \tilde{D}^\dagger_{rig}(-) :\Rep_{L}(G_L) \to \mathfrak{M}(\tilde{\mathcal{R}}_L)$ is exact. Moreover, we have a comparison isomorphism
     \begin{equation}\label{f:comptildeDdaggerrig}
      \tilde{\mathcal{R}}\otimes_{ \tilde{\mathcal{R}}_L} \tilde{D}^\dagger_{rig}(V) \xrightarrow{\cong}  \tilde{\mathcal{R}}\otimes_{o_L} V.
     \end{equation}
\end{lemma}

\begin{proof}
 The comparison isomorphism in the proof  of (an analogue of) \cite[Thm.\ 2.13]{KP} implies the comparison isomorphism
\begin{align*}
       \tilde{\mathcal{R}}\otimes_{ \tilde{\mathcal{R}}_L } \tilde{D}^\dagger_{rig}(V)   & \cong  \tilde{\mathcal{R}}\otimes_{o_L} V
\end{align*}
together with the identity $V=(\tilde{\mathcal{R}}\otimes_{ \tilde{\mathcal{R}}_L } \tilde{D}^\dagger_{rig}(V))^{\varphi_L=1}.$ On the other hand the comparison isomorphism \eqref{f:compdagger3} induces by base change an isomorphism
\begin{align*}
  \tilde{\mathcal{R}}\otimes_{ \tilde{\mathbf{B}}^\dagger_L} \tilde{D}^\dagger(V) \xrightarrow{\cong}  \tilde{\mathcal{R}}\otimes_{o_L} V.
\end{align*}
Taking $H_L$-invariants gives the first claim. The exactness of the functor   $ \tilde{D}^\dagger_{rig}(-)$ follows from the exactness of the functor   $ \tilde{D}^\dagger(-)$ by Prop.\  \ref{Dbasechange}.
\end{proof}

Let $R$ be ${\mathbf{B}}_L,$  ${\mathbf{B}}_L^{\dagger}$,  ${\mathcal{R}}_L$, $\tilde{\mathbf{B}}_L,$  $\tilde{\mathbf{B}}_L^{\dagger}$, $\tilde{\mathcal{R}}_L$ and let  correspondingly $R^{int}$ be ${\mathbf{A}}_L,$  ${\mathbf{A}}_L^{\dagger}$,   ${\mathbf{A}}_L^{\dagger}$,  $\tilde{\mathbf{A}}_L,$ $\tilde{\mathbf{A}}_L^{\dagger}$, $\tilde{\mathbf{A}}_L^{\dagger}$.
%An   $\varphi$-module $M$ over $R$ is called \'{e}tale, if $M$ is  a finitely generated free $R$-module $M$ equipped with a  semilinear endomorphism $\varphi_M$ inducing an isomorphism $\varphi_M^{lin}:R\otimes_{R,\varphi_R}M\xrightarrow{\cong} M.$ Such $M$ is called \'{e}tale, if
We denote by $\Phi(R)^{\acute{e}t}$
%(as  subcategory of $\Phi^{\acute{e}t,f}(\tilde{\mathbf{B}}_L^{\dagger} )$  )
the essential image of the base change functor $R\otimes_{R^{int}} -:\Phi^{\acute{e}t,f}(R^{int} )\to \Phi^{\acute{e}t,f}(R)$  (sic!).

 %consisting of such $M$ there exists a $\varphi_L$-module $N$ over $\tilde{\mathbf{A}}_L^\dagger$, such that  $\tilde{\mathcal{R}}_L\otimes_{\tilde{\mathbf{A}}_L^\dagger} N\cong M$.

\begin{proposition}\label{prop:basechangeRobba}
Base change induces an equivalence of categories
\[
\Phi(\tilde{\mathbf{B}}_L^{\dagger} )^{\acute{e}t} \leftrightarrow \Phi(\tilde{\mathcal{R}}_L^{})^{\acute{e}t}
\]
 and  an isomorphism of Yoneda extension groups
\[
\mathrm{Ext}^1_{\Phi(\tilde{\bf B}_L^{\dagger} )}(\tilde{\bf B}_L^{\dagger},M) \cong
\mathrm{Ext}^1_{\Phi(\tilde{\mathcal{R}}_L )}(\tilde{\mathcal{R}}_L, \tilde{\mathcal{R}}_L \otimes_{\tilde{\bf B}_L^{\dagger}}  M )
\]
for all $M \in \Phi(\tilde{\bf B}_L^{\dagger} )^{\acute{e}t}.$
\end{proposition}

\begin{proof}
The first claim is an analogue of \cite[Thm.\ 8.5.6]{KLI}. The second claim follows as in the proof of Prop.\  \eqref{prop:basechange} using  the fact that by Lemma 8.6.3 in loc.\ cit.\ any extension of \'{e}tale $\varphi$-modules over $\tilde{R}_L$ is again \'{e}tale. Note that $\tilde{\mathcal{R}}_L/\tilde{\mathbf{B}}^\dagger_L$ is a  faithfully flat ring extension, $ \tilde{\mathbf{B}}^\dagger_L$ being a field.
\end{proof}

\begin{corollary}
If $V$ belongs to
  $Rep_{L}(G_L)$, the following complex concentrated  in degrees $0$ and $1$ is acyclic
     \begin{align}\label{f:phidaggertilde}
        &  \xymatrix@C=0.5cm{
          0 \ar[r] & \tilde{D}^\dagger_{rig}(V)/\tilde{D}^\dagger(V) \ar[rr]^{\varphi-1} && \tilde{D}^\dagger_{rig}(V)/\tilde{D}^\dagger(V) \ar[r] & 0 .}
     \end{align}
     In particular, we have that the $n$th cohomology groups of the complex concentrated  in degrees $0$ and $1$
     \[ \xymatrix@C=0.5cm{
          0 \ar[r] & \tilde{D}^\dagger_{rig}(V) \ar[rr]^{\varphi-1} && \tilde{D}^\dagger_{rig}(V) \ar[r] & 0 }\] are isomorphic to $H^n(H_L,V)$ for $n\geq 0.$
\end{corollary}

\begin{proof}
Compare with \cite[Thm.\ 8.6.4]{KLI}  and its proof (Note that the authors meant to cite Thm.\  8.5.12 (taking
c=0, d=1) instead of Thm.\  6.2.9 - a reference which just does not exist within that book). Using the interpretation of the $H^i_{\varphi}$ as $\Hom$- and $\Ext^1$-groups, respectively, the assertion is immediate from Prop.\  \ref{prop:basechangeRobba}. The last statement now follows from Corollary \ref{cor:phitildedagger}.
\end{proof}

\begin{proposition}
 Base extension  gives rise to an equivalence of categories
 \[\mathfrak{M}^{\acute{e}t}({\mathbf{B}}_L^\dagger) \leftrightarrow \mathfrak{M}^{\acute{e}t}({\cR}_L) .\]
\end{proposition}

\begin{proof}
\cite[Prop.\ 1.6]{FX}.
\end{proof}

\begin{lemma}\phantomsection\label{lem:stronghypKed}
\begin{enumerate}
\item ${\bf B}_L^{\dagger}\subseteq {\cR}_{L}$ are B\'{e}zout domains and the strong hypothesis in the sense of \cite[Hypothesis 1.4.1]{Ked08} holds, i.e., for any $n\times n$ matrix $A$ over ${\bf A}_L^{\dagger}$  the map $({\cR}_{L}/{\bf B}_L^{\dagger})^n \xrightarrow{1-A\varphi_L}  ({\cR}_{L}/{\bf B}_L^{\dagger})^n$ is bijective.
%\item ${\cR}_{L}^\times\subseteq {\bf B}_L^{\dagger}.$
\end{enumerate}
\end{lemma}

\begin{proof}
  \cite[Prop.\ 1.2.6]{Ked08}.
\end{proof}

\begin{proposition}
If $V$ belongs to
  $Rep_{L}^\dagger(G_L)$, the following complex concentrated  in degrees $0$ and $1$ is acyclic
     \begin{align}\label{f:phidagger}
        &  \xymatrix@C=0.5cm{
          0 \ar[r] & {D}^\dagger_{rig}(V)/{D}^\dagger(V) \ar[rr]^{\varphi-1} && {D}^\dagger_{rig}(V)/{D}^\dagger(V) \ar[r] & 0 ,}
     \end{align}
     where ${D}^\dagger_{rig}(V):={\cR}_L\otimes_{{\mathbf{B}}_L^\dagger}{D}^\dagger(V).$
     In particular,  the complexes concentrated  in degrees $0$ and $1$
     \[ \xymatrix@C=0.5cm{
          0 \ar[r] & {D}^\dagger_{rig}(V) \ar[rr]^{\varphi-1} && {D}^\dagger_{rig}(V) \ar[r] & 0 } \mbox{ and } \xymatrix@C=0.5cm{
          0 \ar[r] & {D}^\dagger(V) \ar[rr]^{\varphi-1} && {D}^\dagger(V) \ar[r] & 0 }\]   have   the same cohomology groups of  for $n\geq 0$.
\end{proposition}

\begin{proof}
 This follows from the strong hypothesis in Lemma \ref{lem:stronghypKed}   as the Frobenius endomorphism on $M\in \mathfrak{M}^{\acute{e}t}({\mathbf{B}}_L^\dagger)$ is of the form $A\varphi_L$ by definition.
 %because there exists a Frobenius and Galois stable lattice $N$ over ${\mathbf{A}}^\dagger_L$  such that $M={\mathbf{B}}^\dagger_L \otimes_{{\mathbf{A}}^\dagger_L } N.$
\end{proof}

\begin{lemma}
Base change induces fully faithful embeddings $\Phi(\mathbf{A}_L^\dagger)^{\acute{e}t}\subseteq \Phi(\mathbf{A}_L)^{\acute{e}t}$    and $\Phi(\mathbf{B}_L^\dagger)^{\acute{e}t}\subseteq \Phi(\mathbf{B}_L)^{\acute{e}t}$.
\end{lemma}

\begin{proof}
As in the proof of Prop.\  \ref{prop:basechange} this reduces to checking that $\left(\mathbf{A}_L\otimes_{\mathbf{A}_L^\dagger}M\right)^{\varphi=\id}\subseteq M.$ By that proposition we know that
\[\left(\mathbf{A}_L\otimes_{\mathbf{A}_L^\dagger}M\right)^{\varphi=\id}\subseteq \left(\tilde{\mathbf{A}}_L\otimes_{\mathbf{A}_L^\dagger}M\right)^{\varphi=\id} \subseteq \tilde{\mathbf{A}}_L^\dagger\otimes_{\mathbf{A}_L^\dagger} M.\] Since $\mathbf{A}_L\cap \tilde{\mathbf{A}}_L^\dagger=\mathbf{A}_L^\dagger$ within $\tilde{\mathbf{A}}_L$ by definition, the claim follows for the integral version, whence also for the other one my tensoring the integral embedding with $L$ over $o_L.$
\end{proof}

\begin{remark} Note that $H^0_\dagger(H_L,V)=H^0(H_L,V)$ and $H^1_\dagger(H_L,V)\subseteq H^1(H_L,V)$. For the latter relation  use the previous lemma, which implies that an extension which splits after base change already splits itself, together with Corollary \ref{cor:phi} and Remark \ref{rem:Yoneda}. In general the inclusion for $H^1$ is strict as follows indirectly from \cite{FX}. Indeed, otherwise  the complex  \begin{align}\label{f:phidaggertildefalse}
        &  \xymatrix@C=0.5cm{
          0 \ar[r] & {D}(V)/{D}^\dagger(V) \ar[rr]^{\varphi-1} && {D}(V)/{D}^\dagger(V) \ar[r] & 0 ,}
     \end{align}
would be always acyclic, which would imply by the same observation as in Prop.\  \ref{prop:herr} below together with \cite[Thm.\  5.2.10(ii)]{SV23} that $H^1_\dagger(G_L,V)= H^1(G_L,V)$ in contrast to Remark 5.2.13 in (loc.\ cit.).
\end{remark}

\section{The  web of eqivalences}\label{subsec:picture}

We summarize the various equivalences of  categories, for which we only sketch proofs or indicate analogue results whose proofs can be transferred to our setting.

\begin{proposition}\label{prop:equivalenceintegral}
The following categories are equivalent:
\begin{enumerate}
\item $Rep_{o_L}(G_L)$,
\item  $\mathfrak{M}^{\acute{e}t}({\mathbf{A}}_L) ,$
\item  $\mathfrak{M}^{\acute{e}t}(\tilde{\mathbf{A}}_L) $ and
\item $\mathfrak{M}^{\acute{e}t}(\tilde{\bf A}_L^{\dagger} ).$
\end{enumerate}
The equivalences from $(ii)$ and $(iv)$ to $(iii)$ are induced by base change.
\end{proposition}

\begin{proof}
This can be proved in the same way as in \cite[Thm.\ 2.3.5]{KedNew}, although it seems to be only a sketch. Another way is to check that the very detailed proof for the equivalence between (i) and (ii) in \cite{GAL} almost literally carries over to a proof for the equivalence between (i) and (iii).  Alternatively, this is a consequence of Prop.\  \ref{prop:weakdecomp} by \cite[Thm.\ 5.4.6]{KLII}. See also \cite{Kl}. For the equivalence between (iii) and (iv) consider the $2$-commutative diagram
\[\xymatrix{
  \mathfrak{M}^{\acute{e}t}(\tilde{\bf A}_L^{\dagger} )   \ar[rr]^{\text{faithfully flat base change}} &   & \mathfrak{M}^{\acute{e}t}(\tilde{\bf A}_L^{} ) \ar[dl]^{ } \\
   & \Rep_{o_L}(G_L) \ar[ru]^{ } \ar[lu] &    , }\]
which is induced by the isomorphism \eqref{f:basechangedagger2}  and immediately implies (essential) surjectivity on objects and morphisms while the faithfulness follows from faithfully flat base change.
\end{proof}

\begin{corollary}\label{cor:equivalence}
The following categories are equivalent:
\begin{enumerate}
\item $Rep_{L}(G_L)$,
\item  $\mathfrak{M}^{\acute{e}t}({\mathbf{B}}_L) $,
\item  $\mathfrak{M}^{\acute{e}t}(\tilde{\mathbf{B}}_L) $ and
\item  $\mathfrak{M}^{\acute{e}t}(\tilde{\mathbf{B}}_L^\dagger) $.
\end{enumerate}
The equivalences from $(ii)$  and $(iv)$ to $(iii)$ are induced by base change.
\end{corollary}

\begin{proof}
This follows from Propositions \ref{prop:basechangephiGamma} and \ref{prop:equivalenceintegral} by inverting $\pi_L.$
\end{proof}

\begin{proposition}\label{prop:equivalence}
The categories in Corollary \ref{cor:equivalence} are - via base change from (iv) - also equivalent to
\begin{enumerate}
\item[{\rm (v)}] $ \mathfrak{M}^{\acute{e}t}(\tilde{\mathcal{R}}_L) .$
\end{enumerate}
\end{proposition}

\begin{proof}
By definition base change is essentially surjective and it    is well-defined - regarding the continuity of the $\Gamma_L$-action - by  Lemma \ref{lem:LF} and Lemma \ref{lem:contGammaAction}. Since for   \'{e}tale $\varphi_L$-modules we know fully faithfulness already, taking $\Gamma_L$-invariants gives fully faithfulness for $(\varphi_L,\Gamma_L)$-modules, too. \footnote{Regarding $\varphi_L$-modules cf.\ \cite[the equivalence between (e) and (f) of Thm.\ 8.5.6]{KLI}, see also  Thm.\ 8.5.3 in (loc.\ cit.), the equivalence (d) to (e).}
\end{proof}

Altogether we may visualize the relations between the various categories by the following diagram:
%\begin{equation*}
%   \xymatrix{
%     \mathfrak{M}^{\acute{e}t}(\tilde{\mathcal{R}}_L)   & \mathfrak{M}^{\acute{e}t}(\tilde{\mathbf{B}}_L^\dagger)  \ar@{=>}[l]^{ }\ar@{=>}[r]^{ } & \mathfrak{M}^{\acute{e}t}(\tilde{\mathbf{B}}_L)  \ar@{<=>}[rd]^{\tilde{V}}_{\tilde{D}} \\
%     \mathfrak{M}^{\acute{e}t}({\mathcal{R}}_L)\ar[u]_{ }  &  \mathfrak{M}^{\acute{e}t}({\mathbf{B}}_L^\dagger)  \ar[u]_{ }\ar@{=>}[l]^{ }\ar[r]^{ } & \mathfrak{M}^{\acute{e}t}({\mathbf{B}}_L) \ar@{=>}[u]^{ } &  \ar@{<=>}[l]^{D}_{V}  \Rep_{L}(G_L)}
%\end{equation*}
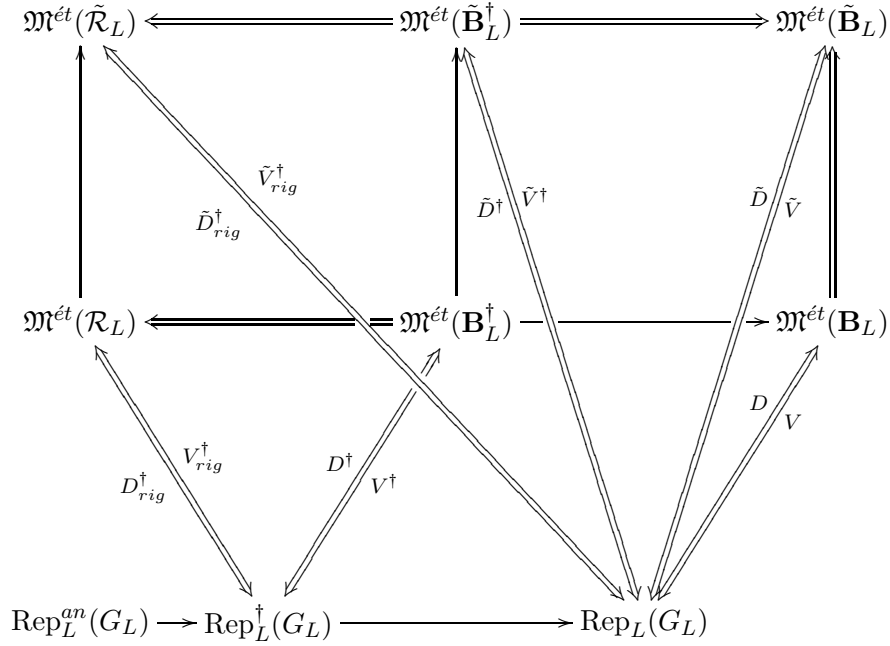
\begin{equation*}\label{diag1}
\begin{xy}
%vertices
(25,-20)*+{ \Rep_{L}^\dagger(G_L)  }="c1";           (75,-20)*+{ \Rep_{L}(G_L)}="c2";%
(0,-20)*+{ \Rep_{L}^{an}(G_L)  }="c3";
(0,20)*+{\mathfrak{M}^{\acute{e}t}({\mathcal{R}}_L)}="b1";   (50,20)*+{\mathfrak{M}^{\acute{e}t}({\mathbf{B}}_L^\dagger)}="b2";    (100,20)*+{\mathfrak{M}^{\acute{e}t}({\mathbf{B}}_L) }="b3"; %
(0,60)*+{\mathfrak{M}^{\acute{e}t}(\tilde{\mathcal{R}}_L) }="a1";   (50,60)*+{\mathfrak{M}^{\acute{e}t}(\tilde{\mathbf{B}}_L^\dagger) }="a2";    (100,60)*+{ \mathfrak{M}^{\acute{e}t}(\tilde{\mathbf{B}}_L)  }="a3";
%arrows
{\ar@{->} "c3"; "c1"};{\ar@{->} "c1"; "c2"};{  \ar@{<=>}^(0.7){D}_(0.7){V} "c2"; "b3" }; {\ar@{<=>}^(0.3){\tilde{V}}_(0.3){\tilde{D}}  "a3";"c2"};
{  \ar@{<=>}^{D^\dagger}_{V^\dagger}|(.8){\hole} "c1"; "b2" };
{  \ar@{<=>}^{D^\dagger_{rig}}_{V^\dagger_{rig}} "c1"; "b1" };
{\ar@{<=>}^(0.3){\tilde{V}^\dagger}_(0.3){\tilde{D}^\dagger}  "a2";"c2"}; {\ar@{<=>}^(0.3){\tilde{V}^\dagger_{rig}}_(0.3){\tilde{D}^\dagger_{rig}}  "a1";"c2"}; {\ar@{=>}"a2"; "a1"};{\ar@{=>}"a2"; "a3"};
{\ar@{->} "b1"; "a1"};{\ar@{->} "b2"; "a2"};{\ar@{=>}|(.25){\hole} "b2"; "b1"};{\ar@{->}|(.25){\hole}|(.75){\hole}"b2"; "b3"};{\ar@{=>}"b3"; "a3"};
\end{xy}
\end{equation*}

Here all arrows represent functors which are fully faithful, i.e., embeddings of categories. Arrows without label denote base change functors. Under them the functors $D, \tilde{D}, {D}^\dagger, \tilde{D}^\dagger, {D}^\dagger_{rig}$, and $\tilde{D}^\dagger_{rig}$ are compatible. The arrows ${=>} $ represent equivalences of categories, while the arrows ${-\hspace{-0.15cm}>}$ represent embeddings which are not essentially surjective in general. We recall that the quasi-inverse functors are given as follows
\begin{align*}
V(M)=&(\mathbf{B}\otimes_{\mathbf{B}_L}M)^{\varphi=1},\;\;  \tilde{V}(M)=(\tilde{\mathbf{B}}\otimes_{\tilde{\mathbf{B}}_L}M)^{\varphi=1} ,\;\; V^\dagger(M)=(\mathbf{B}^\dagger\otimes_{\mathbf{B}_L^\dagger}M)^{\varphi=1},\\
\tilde{V}^\dagger(M)=&(\tilde{\mathbf{B}}^\dagger\otimes_{\tilde{\mathbf{B}}_L^\dagger}M)^{\varphi=1},\;\;
\tilde{V}^\dagger_{rig}(M)= (\tilde{\mathcal{R}}\otimes_{ \tilde{\mathcal{R}}_L } M)^{\varphi=1}\mbox{ and }
{V}^\dagger_{rig}(M)= (\tilde{\mathcal{R}}\otimes_{ {\mathcal{R}}_L } M)^{\varphi=1}.
\end{align*}\footnote{By \cite[Prop.\ 1.5 (a)]{FX} the third formula holds while by (c) there is an equivalence of categories.}
\footnote{For the fourth formula compare with the proof of Propositon \ref{prop:basechange} omitting the index $L$ in $\tilde{\mathbf{A}}^\dagger_L$, etc. to conclude that $(\tilde{\mathbf{B}}\otimes_{\mathbf{B}_L^\dagger}M)^{\varphi=1}=(\tilde{\mathbf{B}}^\dagger\otimes_{\mathbf{B}_L^\dagger}M)^{\varphi=1}$.}
\footnote{Since $V^\dagger(M_0)\subseteq {V}^\dagger_{rig}(M)\subseteq\tilde{V}^\dagger_{rig}(\tilde{\mathcal{R}}_L\otimes_{\mathcal{R}_L}M) $ for some model $M_0$ over $\mathbf{B}^\dagger_L$ of $M$ we obtain the last formula.}

\section{Cohomology: Herr complexes}\label{sec:cupprod}
The aim of this section is to compare the Herr complexes of the various $(\varphi_L,\Gamma_L)$-modules attached to a given Galois representation.

We fix some open subgroup $U\subseteq \Gamma_L$ and let $L'=L_\infty^U.$

Let $M_0$ be a complete linearly topologised $o_L$-module with continuous  $U$-action and with continuous $U$-equivariant endomorphism $f$.  We define \[\mathcal{T} := \mathcal{T}_{f,U}(M_0):=\mathrm{cone}\left( \cC^\bullet(U,M_0)\xrightarrow{(f)_*-1} \cC^\bullet(U,M_0)\right)[-1]\]   the mapping fibre of  $\cC^\bullet(U,f-1)$. The importance of this generalized Herr-complex is given by the fact that it computes Galois cohomology when applied to $M_0=D(V)$ and $f=\varphi_{D(V)}:$

\begin{theorem}\label{thm:Herr}
Let $V$ be in $\Rep_{L}(G_L)$
 For $  D(V)$ the corresponding $(\varphi_L,\Gamma_L)$-module over $\mathbf{B}_L$ we have canonical isomorphisms
\begin{equation}\label{f:cohLphiB}
  h^* = h^*_{U,V}: H^*(L',V) \xrightarrow{\ \cong\ }   h^*(\mathcal{T}_{\varphi,U}(D(V)) )
\end{equation}
which are functorial in $V$ and compatible with restriction and corestriction.
\end{theorem}
  \begin{proof}
To this aim let $T$ be a $G_L$-stable lattice of $V$.
In \cite[Thm.\ 5.1.11.]{Ku}, \cite[Thm.\ 5.1.11.]{KV}  it is shown that the cohomology groups of $ \mathcal{T}_{\varphi,U}({D}(T))$ are canonically isomorphic to $H^i(L',T)$ for all $i\geq 0,$ whence the cohomology groups of $ \mathcal{T}_{\varphi,U}({D}(T))[\frac{1}{\pi_L}]$ are canonically isomorphic to $H^i(L',V)$ for all $i\geq 0$.
\end{proof}
  Note that we obtain a decomposition $U\cong \Delta\times U'$ with
  a subgroup $U'\cong\mathbb{Z}_p^d$ of $U$ and  $\Delta$  the torsion subgroup of $U$.
We now fix topological generators $\gamma_1,\ldots \gamma_d$  of $U'$ and we set $\Lambda:= \Lambda(U').$ By \cite[Thm.\ II.2.2.6]{La} the $U'$-actions extends to continuous $\Lambda$-action and one has $\mathrm{Hom}_{\Lambda, cts}(\Lambda,M_0)=\mathrm{Hom}_{\Lambda}(\Lambda,M_0).$
Consider the (homological) complexes $K_\bullet(\gamma_i):=[\Lambda\xrightarrow{\gamma_i-1}\Lambda]$ concentrated in degrees $1$ and $0$ and define the Koszul complexes
\begin{align*}
K_\bullet:=& K_\bullet^{U'}:= K_\bullet(\gamma):={\bigotimes_{{\substack{\phantom{i=1}{\Lambda}\\i=1}}}^d}K_\bullet(\gamma_i) \;\;\mbox{ and }\\
K^\bullet(M_0):=&K^\bullet_{U'}(M_0):=\Hom_\Lambda^{\bullet} (K_\bullet,M_0)\cong \Hom_\Lambda^{\bullet }(K_\bullet,\Lambda)\otimes_\Lambda M_0=K^\bullet(\Lambda)\otimes_\Lambda M_0.
%K_\bullet(M_0):=&K_\bullet\otimes_\Lambda M_0 \mbox{ (homological complex)},\\
%K_\bullet(M_0)^\bullet:=&(K_\bullet\otimes_\Lambda M_0)^\bullet \mbox{ (the associated cohomological complex)}.
\end{align*}
Following \cite[\S 4.2]{CoNi} and \cite[(169)]{SV23} we obtain     a quasi-isomorphism
\begin{equation}\label{f:KoszulCtsCoh}
  K^\bullet_{U'}(M_0)\xrightarrow{\simeq} \cC^\bullet(U',M_0)
\end{equation}
 inducing the quasi-isomorphism
\begin{equation}\label{f:KoszulCtsCohPhi}
  K_{f,U'}(M_0)\xrightarrow{\simeq} \mathcal{T}_{f,U'}(M_0),
\end{equation}
where we denote by $K_{f,U'}(M_0):=\mathrm{cone}\left( K^\bullet(M_0)\xrightarrow{f-\id}K^\bullet(M_0) \right)[-1]$ the mapping fibre of  $K^\bullet(f)$.
%\begin{lemma}\label{lem:HS}
% There is a canonical quasi-isomorphism \[\mathcal{T}_{f,U}(M_0)[\frac{1}{\pi_L}]\simeq K^\bullet_{f,U'}(M_0[\frac{1}{\pi_L}]^\Delta).\] If $M_0$ is an $L$-vector space, the inversion of $\pi_L$ can be omitted on both sides.
%\end{lemma}
More generally, by  \cite[Lem.\ A.0.1]{SV23} we obtain a canonical quasi-isomorphism
  \begin{equation}\label{f:KoszulCtsCohPhiDelta}
  K_{f,U'}(M^\Delta)\xrightarrow{\simeq} \mathcal{T}_{f,U}(M),
\end{equation}
i.e., by Theorem \ref{thm:Herr} we also have   canonical isomorphisms
\begin{equation}\label{f:cohLphiB2}
  h^* = h^*_{U,V}: H^*(L',V) \xrightarrow{\ \cong\ }   h^*(K_{f,U'}(D(V)^\Delta) ).
\end{equation}

The next proposition extends this result to $\tilde{D}(V), \tilde{D}^\dagger(V)$ and $ \tilde{D}^\dagger_{rig}(V)$ instead of ${D}(V)$.

\begin{proposition}\label{prop:herr}
If $V$ belongs to $Rep_{L}(G_L)$,   the canonical inclusions of Herr complexes \begin{align*}\label{}
 K^\bullet_{\varphi,U'}(\tilde{D}^\dagger(V)^\Delta) &\subseteq  K^\bullet_{\varphi,U'}(\tilde{D}^\dagger_{rig}(V)^\Delta) ,\\
 K^\bullet_{\varphi,U'}(\tilde{D}^\dagger(V)^\Delta) &\subseteq  K^\bullet_{\varphi,U'}(\tilde{D}(V)^\Delta) \mbox{ and }\\
K^\bullet_{\varphi,U'}({D}(V)^\Delta) &\subseteq  K^\bullet_{\varphi,U'}(\tilde{D}(V)^\Delta)\label{f:DtildeD}
\end{align*}     are quasi-isomorphisms and their cohomology groups are canonically isomorphic to $H^i(L',V)$ for all $i\geq 0.$
\end{proposition}

\begin{proof}
 Forming Koszul complexes with regard to $U'$ we obtain  the following diagram of (double) complexes with exact columns
\[\xymatrix{
0\ar[d] & 0\ar[d]\\
  K^\bullet({D}(V)^\Delta) \ar[d]_{} \ar[r]^{\varphi-1} &  K^\bullet({D}(V)^\Delta)\ar[d]^{} \\
   K^\bullet( \tilde{D}(V)^\Delta)\ar[d]_{} \ar[r]^{\varphi-1} &  K^\bullet(\tilde{D}(V)^\Delta) \ar[d]^{} \\
  K^\bullet((\tilde{D}(V)/{D}(V))^\Delta)\ar[d] \ar[r]^{\varphi-1}_{\cong} & K^\bullet( (\tilde{D}(V)/{D}(V))^\Delta)\ar[d] \\
  0& 0  }\]
in which the bottom line is an isomorphism of complexes by \ref{cor:phi}, as the action of $\Delta$ commutes with $\varphi$. Hence,
going over to total complexes gives  an exact sequence
\[ 0\to K_{\varphi,U}^\bullet({D}(V)^\Delta) \to K_{\varphi,U}^\bullet(\tilde{D}(V)^\Delta)\to K_{\varphi,U}^\bullet((\tilde{D}(V)/{D}(V))^\Delta)\to 0, \]
in which $K_{\varphi,U}^\bullet((\tilde{D}(V)/{D}(V))^\Delta) $ is acyclic. Thus we have shown the statement regarding the last inclusion. The other two cases are dealt with similarly, now using  \eqref{f:phidaggertilde} and \ref{cor:phitildedagger} combined with \eqref{f:phitilde}. It follows in particular that all six Koszul complexes in the statement are quasi-isomorphic. Therefore   the second part of the assertion  follows from \eqref{f:cohLphiB2}.
\end{proof}

In accordance with diagram at the end of subsection \ref{subsec:picture}  we may visualize the relations between the various Herr complexes by the following diagram:
%\begin{equation*}
%   \xymatrix{
%     \mathfrak{M}^{\acute{e}t}(\tilde{\mathcal{R}}_L)   & \mathfrak{M}^{\acute{e}t}(\tilde{\mathbf{B}}_L^\dagger)  \ar@{=>}[l]^{ }\ar@{=>}[r]^{ } & \mathfrak{M}^{\acute{e}t}(\tilde{\mathbf{B}}_L)  \ar@{<=>}[rd]^{\tilde{V}}_{\tilde{D}} \\
%     \mathfrak{M}^{\acute{e}t}({\mathcal{R}}_L)\ar[u]_{ }  &  \mathfrak{M}^{\acute{e}t}({\mathbf{B}}_L^\dagger)  \ar[u]_{ }\ar@{=>}[l]^{ }\ar[r]^{ } & \mathfrak{M}^{\acute{e}t}({\mathbf{B}}_L) \ar@{=>}[u]^{ } &  \ar@{<=>}[l]^{D}_{V}  \Rep_{L}(G_L)}
%\end{equation*}
\begin{equation*}\label{diag2}
\begin{xy}
%vertices
(25,-20)*+{   }="c1";           (75,-20)*+{ \cC^\bullet(G_{L'},V)}="c2";
(0,-20)*+{    }="c3";
(0,20)*+{K^\bullet_{\varphi,U'}({D}^\dagger_{rig}(V)^\Delta)}="b1";   (50,20)*+{K^\bullet_{\varphi,U'}({D}^\dagger(V)^\Delta)}="b2";    (100,20)*+{K^\bullet_{\varphi,U'}({D}(V)^\Delta) }="b3"; %
(0,60)*+{K^\bullet_{\varphi,U'}(\tilde{D}^\dagger_{rig}(V)^\Delta) }="a1";   (50,60)*+{K^\bullet_{\varphi,U'}(\tilde{D}^\dagger(V)^\Delta) }="a2";    (100,60)*+{ K^\bullet_{\varphi,U'}(\tilde{D}(V)^\Delta)  }="a3";
%arrows
%{\ar@{->} "c3"; "c1"};{\ar@{->} "c1"; "c2"};
{  \ar@{<-->}  "c2"; "b3" };
%{\ar@{<=>}^(0.3){\tilde{V}}_(0.3){\tilde{D}}  "a3";"c2"};
%{  \ar@{<=>}^{D^\dagger}_{V^\dagger}|(.8){\hole} "c1"; "b2" };
%{  \ar@{-->} "c1"; "b1" };
%{\ar@{<=>}^(0.3){\tilde{V}^\dagger}_(0.3){\tilde{D}^\dagger}  "a2";"c2"}; {\ar@{<=>}^(0.3){\tilde{V}^\dagger_{rig}}_(0.3){\tilde{D}^\dagger_{rig}}  "a1";"c2"};
{\ar@{=>}"a2"; "a1"};{\ar@{=>}"a2"; "a3"};
{\ar@{->}^{} "b1"; "a1"};{\ar@{->} "b2"; "a2"};{\ar@{=>}|(.25){\hole} "b2"; "b1"};{\ar@{->}|(.25){\hole}|(.75){\hole}"b2"; "b3"};{\ar@{=>}"b3"; "a3"};
\end{xy}
\end{equation*}

Here all arrows represent injective maps of complexes, among which the arrows ${=>} $ represent quasi-isomorphisms, while the arrows ${-\hspace{-0.15cm}>}$ need not induce isomorphisms on cohomology, in general. The interrupted arrow ${-\;-\hspace{-0.15cm}>}$  means a map in the derived category while  ${<-\;-\hspace{-0.15cm}>}$ means a quasi-isomorphism in the derived category. By \cite[Lem.\ A.0.1]{SV23} we have a analogous diagram for $\cT_{\varphi,U}(?(V))$ with $?\in\{D, \tilde{D}, {D}^\dagger, \tilde{D}^\dagger, {D}^\dagger_{rig}, \tilde{D}^\dagger_{rig}\}.$

\begin{remark}
The image of \[h^i(\cT_{\varphi,U}(D^\dagger_{rig}(V)))\cong h^i(K^\bullet_{\varphi,U'}({D}^\dagger_{rig}(V)^\Delta))\cong h^i(K^\bullet_{\varphi,U'}({D}^\dagger(V)^\Delta))\cong h^i(\cT_{\varphi,U}(D^\dagger(V)))\] in $H^i(L',V)$ is independent of the composite ($=$ path) in above diagram.
\end{remark}

\section{Weakly decompleting towers}\label{sec:weak}

 Kedlaya and Liu's developed in \cite[\S 5]{KLII} the concept of perfectoid towers and studied their properties in an axiomatic way. The aim of this section is to show that the Lubin-Tate extensions considered in this article form  a {\it weakly decompleting}, but not a {\it decompleting} tower, properties which we will recall or refer to in the  course of this section. Moreover, we have to show that the  {\it axiomatic} period rings coincide with those introduced earlier.\\

In the sense of Def.\ 5.1.1 in (loc.\ cit.) the sequence $\Psi=(\Psi_n: (L_n, o_{L_n})\to (L_{n+1},o_{L_{n+1}}))_{n=0}^\infty$ forms a finite \'{e}tale tower over $(L,o_L)$ or $X:=\mathrm{Spa}(L,o_L)$, which is perfectoid as $\hat{L}_\infty$ is by \cite[Prop.\ 1.4.12]{GAL}.\footnote{ In the notation of \cite{KLII}:
$E=L,$ $\varpi=\pi_L$, $h=r$, $k:=o_L/(\pi_L)=\mathbb{F}_q,$  i.e. $q=p^r$.  $A_{\Psi,n}:=L_n,$  $A_{\Psi,n}^+:=o_{L_n},$ $X:=\mathrm{Spa}(L,o_L)$ with the obvious transition maps which are finite \'{e}tale.

$(A_\Psi, A_{\Psi}^+):=\varinjlim_n (A_{\Psi,n},A_{\Psi,n}^+)= (L_\infty, o_{L_\infty})$

$(\tilde{A}_\Psi,\tilde{A}_\Psi^+):=(A_\Psi, A_{\Psi}^+)^{\wedge \pi_L-adic}=(\hat{L}_\infty, o_{\hat{L}_\infty})$
   }

Therefore we can use the perfectoid correspondence \cite[Thm.\ 3.3.8]{KLII} to associate with $ (\hat{L}_\infty, o_{\hat{L}_\infty})$ the pair
\[(\tilde{R}_\Psi,\tilde{R}_\Psi^+):=(\hat{L}_\infty^\flat, o_{\hat{L}_\infty}^\flat).\] Now we recall the variety of period rings, which Kedlaya and Liu attach to the tower, in our notation,  starting with\\

{\bf Perfect period rings:}\Footnote{ Eigentlich sheaves, wie ist Folgendes zu verstehen: \\  Convention 5.0.2. Throughout \S5, we work over an affinoid space rather than an arbitrary
preadic space; we thus freely use the {\bf A;B;C} suite of notations to refer to period rings rather
than period sheaves on $X_{pro\acute{e}t}$.\\

Welches $X$ genau und welche Schnitte werden hier nun ausgew\"{a}hlt?}
\begin{align*}
  \tilde{\bf A}_\Psi & :=\tilde{\bf A}_L=W(\hat{L}_\infty^\flat)_L, \\
 \tilde{\bf A}_\Psi^+ & :=W(o_{\hat{L}_\infty}^\flat)_L\subseteq
\tilde{\bf A}_\Psi^{\dagger,r}:=\tilde{\bf A}_L^{\dagger,r}=\{x=\sum_{i\geq 0} \pi_L^i[x_i]\in W({\hat{L}_\infty}^\flat)_L| \mbox{ }  |\pi_L^{i}\|x_i|_\flat^r\xrightarrow{i\to\infty} 0\}, \\
 \tilde{\bf A}_\Psi^{\dagger}& :=\bigcup_{r>0} \tilde{\bf A}_\Psi^{\dagger,r}=\tilde{\bf A}_L^{\dagger}
\end{align*}
%$\tilde{\bf A}_\Psi:=\tilde{\bf A}_L=W(\hat{L}_\infty^\flat)_L,$ $\tilde{\bf A}_\Psi^+:=W(o_{\hat{L}_\infty}^\flat)_L\subseteq
%\tilde{\bf A}_\Psi^{\dagger,r}:=\tilde{\bf A}_L^{\dagger,r}=\{x=\sum_{i\geq 0} \pi_L^i[x_i]\in W({\hat{L}_\infty}^\flat)_L| \mbox{ for $i$ to $\infty$ :}  |\pi_L^{i}\|x_i|_\flat^r\rightarrow 0\},$ $\tilde{\bf A}_\Psi^{\dagger}:=\bigcup_{r>0} \tilde{\bf A}_\Psi^{\dagger,r}=\tilde{\bf A}_L^{\dagger}$.\\
%\footnote{
%\com{\color{blue} Meiner Meinung nach gilt
%\begin{equation*}
%  \tilde{\bf A}_\Psi^{\dagger,r} = \tilde{\bf A}^{\dagger,\frac{q-1}{qr}}_L \ ,
%\end{equation*}
%wobei die rechte Seite die Ringe sind, die wir in Heidelberg betrachtet haben, und ich den Normalisierungen Deines
%Studenten folge mit $|\pi_L| = q^{-1}$. Deswegen bin ich daf\"ur, dass wir nur ein System von Normen benutzen. Im \"Ubrigen scheint mir auch zu gelten:
%\begin{equation*}
%  \tilde{\bf A}_\Psi^{\dagger} = \bigcup_{n > 0} \varphi_L^n(\tilde{\bf A}_\Psi^{\dagger,1}) \ .
%\end{equation*}
%}}

 {\bf Imperfect period rings:} \\

 To introduce these we first recall the map $\Theta: W(o_{\mathbb{C}_p}^\flat)_L\to o_{\mathbb{C}_p},\; \sum_{i\geq 0} \pi_L^i [x_i]\mapsto \sum \pi_L^ix_i^\sharp,$ which extends to a map $\Theta: \tilde{\bf A}_\Psi^{\dagger,s}\to {\mathbb{C}_p}$ for all $s\geq1;$ for arbitrary $r>0$ and $n\geq -\log_q r$ the composite  $\tilde{\bf A}_\Psi^{\dagger,r}\xrightarrow{\varphi_L^{-n}} \tilde{\bf A}_\Psi^{\dagger,1}\xrightarrow{\Theta} \mathbb{C}_p$ is well defined and continuous as it is easy to check.   It is a homomorphism of $o_L$-algebras by \cite[Lem.\ 1.4.18]{GAL}.

Following \cite[\S 5]{KLII}  we set ${\bf A}_\Psi^{\dagger,r}:=\{x \in \tilde{\bf A}_\Psi^{\dagger,r}| \Theta(\varphi_q^{-n}(x))\in L_n \mbox{ for all }  n\geq -\log_q r\},$
$ {\bf A}_\Psi^{\dagger}:=\bigcup_{r>0}  {\bf A}_\Psi^{\dagger,r},$ its completion
${\bf A}_\Psi:=({\bf A}_\Psi^{\dagger})^{\wedge \pi_L-adic}$, and residue field
$R_\Psi:= {\bf A}_\Psi/(\pi_L)=({\bf A}_\Psi^{\dagger})/(\pi_L)\subseteq\tilde{R}_\Psi$, $R_\Psi^+:=R_\Psi\cap \tilde{R}_\Psi^+.$

Note that $\omega_{LT}=\{[\iota(t)]\}\in \tilde{\bf A}_\Psi^+:=W(o_{\hat{L}_\infty}^\flat)_L\subseteq\tilde{\bf A}_\Psi^{\dagger,r}$ for all $r>0$ (in the notation of \cite{GAL}).  \cite[Lem.\ 2.1.12]{GAL}  shows
\[\Theta(\varphi_q^{-n}(\omega_{LT}))=\Theta(\{[\varphi_q^{-n}(\omega)]\})=\lim_{i\to\infty} [\pi_L^i]_\varphi(z_{i+n})=z_n\in L_n,\] where $t=(z_n)_{n\geq 1}$ is a fixed generator of the Tate module $T_\pi$ of the formal Lubin-Tate group and $\omega=\iota(t)\in W(o_{\mathbb{C}_p}^\flat )_L$  is the reduction of $\omega_{LT}$ modulo $\pi_L$ satisfying with ${\bf E}_L=k((\omega))$. Therefore $\omega_{LT}$ belongs to ${\bf A}_\Psi^+:={\bf A}_\Psi\cap  \tilde{\bf A}_\Psi^+$. Then it is clear that first ${\bf A}_L^+:=o_L[[\omega_{LT}]] \subseteq  \tilde{{\bf A}}_\Psi^{\dagger}$  and by the continuity of $\Theta\circ \varphi_L^{-n}$ even
${\bf A}_L^+  \subseteq  {\bf A}_\Psi^{\dagger}$ holds. Since  $\omega_{LT}^{-1}\in \tilde{\bf A}_\Psi^{\dagger,\frac{q-1}{q}} $    by  \cite[Lem.\ 3.10]{Ste} (in analogy with \cite[Cor. II.1.5]{ChCo})  and   $\Theta\circ \varphi_L^{-n}$ is a ring homomorphism, it follows that  $\omega_{LT}^{-1}\in {\bf A}_\Psi^{\dagger,\frac{q-1}{q}}$ and    $o_L[[\omega_{LT}]][\frac{1}{\omega_{LT}}]\subseteq  {\bf A}_\Psi^{\dagger}.$
\begin{lemma}
We have $ R_\Psi^+={\bf E}_L^+$ and  $R_\Psi={\bf E}_L.$
\end{lemma}

\begin{proof}
From the above it follows that ${\bf E}_L\subseteq R_\Psi$, whence ${\bf E}_L^{perf}\subseteq R_\Psi^{perf}  \subseteq  \tilde{R}_\Psi=\hat{L}_\infty^\flat$ the latter being perfect. Since $\widehat{{\bf E}_L^{perf}}=\hat{L}_\infty^\flat$ by \cite[Prop.\ 1.4.17]{GAL} we conclude that
\begin{equation}
\label{f:Rperfdense}\mbox{$R^{perf}_\Psi$ is dense in $\tilde{R}_\Psi.$}
\end{equation}
By \cite[Lem.\ 5.2.2]{KLII} have the inclusion
\[R_\Psi^+\subseteq \{x\in\tilde{R}_\Psi| \mbox{$x=(\bar{x}_n)$ with $\bar{x}_n\in o_{L_n}/(z_1)$ for $n>>1$}\}\stackrel{\text{(*)}}{=}{\bf E}_L^+=k[[\omega]]\]
where the   equality (*)  follows from  work of Wintenberger as recalled in \cite[Prop.\ 1.4.29]{GAL}. Since ${\bf E}_L^+\subseteq\tilde{R}^+_\Psi$ by its construction in (loc.\ cit.), we conclude that $ R_\Psi^+={\bf E}_L^+$.

Since each element of $R_\Psi$ is of the form $\frac{a}{\omega^m}$ with $a\in R_\Psi^+$ and $m\geq 0$ by \cite[Lem.\ 1.4.6]{GAL}\footnote{For $\alpha\in R_\Psi $ there exist $m\geq 0$ such that $|\omega^m\alpha|_\flat\leq 1$, i.e., $\omega^m\alpha\in R_\Psi^+.$ }, we conclude that $ R_\Psi={\bf E}_L.$
\end{proof}
 Thus for each $r>0$ such that $\omega_{LT}^{-1}\in  {\bf A}_\Psi^{\dagger,r},$ reduction modulo $\pi_L$ induces a surjection $ {\bf A}_\Psi^{\dagger,r}\twoheadrightarrow R_\Psi$. Recall that $\Psi$ is called weakly decompleting, if
\begin{enumerate}
\item $R^{perf}_\Psi$ is dense in $\tilde{R}_\Psi.$
\item for some $r>0$ we have a strict surjection ${\bf A}_\Psi^{\dagger,r}\twoheadrightarrow R_\Psi$ induced by the reduction modulo $\pi_L$ for the norms $|-|_{r}$ defined by $|x|_{r}:=\sup_i \{|\pi_L^{i}\| x_i|^r_\flat\}$ for $x=\sum_{i\geq 0} \pi_L^i[x_i]$, and $|-|^r_\flat.$
\end{enumerate}
We recall from \cite[Prop.\ 1.4.3.]{FF} or \cite[Prop.\ 5.1.2 (a)]{KLI} that $|-|_{r}$ is multiplicative.
\begin{proposition}\label{prop:weakdecomp}
 The above tower $\Psi$ is weakly decompleting.
\end{proposition}

\begin{proof} Since \eqref{f:Rperfdense} gives (i), only $(ii)$ is missing:
   Since $\omega_{LT}$ has $[\omega]$ in degree zero of its Teichm\"{u}ller series, we may and do choose $r>0$ such that $|\omega_{LT}-[\omega]|_{r}<|\omega|_\flat^r.$ Then $|\omega_{LT}|_{r}=\max\{|\omega_{LT}-[\omega]|_{r}, |\omega|_\flat^{r}\}=|\omega|_\flat^{r}.$ Consider the quotient norm $\|b\|^{(r)}=\inf_{a\in {\bf A}_\Psi^{\dagger,r},a\equiv b \mod \pi_L}|a|_{r}.$

   Now let $b=\sum_{n\geq n_0} a_n\omega^n\in R_\Psi=k((\omega))$ with $a_{n_0}\neq 0$. Lift each $a_n\neq 0$ to $\breve{a}_n\in o_L^\times$ and set $\breve{a}_n=0$ otherwise. Then, for the lift $x:=\sum_{n\geq n_0} \breve{a}_n\omega_{LT}^n$ of $b$ we have by the multiplicativity of $|-|_{r}$ that \[ \|b\|^{(r)}\leq |x|_{r}=(|\omega_{LT}|_{r})^{n_0}=(|\omega|_\flat^{r})^{n_0}=|b|_\flat^r.\] Since,  the other inequality $|b|_\flat^r\leq \|b\|^{(r)} $ giving by continuity is clear, the claim follows.
\end{proof}

\begin{proposition}\label{prop:Lpsi}
${\bf A}_L={\bf A}_\Psi.$
\end{proposition}

\begin{proof}
Both rings have the same reduction modulo $\pi_L.$ And using that the latter element is not a zero-divisor in any of these rings we conclude inductively, that ${\bf A}_L/\pi_L^n{\bf A}_L ={\bf A}_\Psi/\pi_L^n{\bf A}_\Psi$ for all $n.$ Taking projective limits gives the result.
\end{proof}

\begin{proposition}
${\bf A}_L^\dagger={\bf A}_\Psi^\dagger.$
\end{proposition}

\begin{proof}
By \cite[Lem.\ 5.2.10]{KLII} we have that ${\bf A}_\Psi^\dagger=\tilde{\bf A}_L^\dagger\cap \cR_L.$ On the other hand ${\bf A}_L^\dagger=  (\tilde{\bf A}^\dagger \cap {\bf A})^{H_L}=\tilde{\bf A}^\dagger_L \cap {\bf A}$ is contained in $\cR_L$ by Remark \ref{rem:overconvergent}, whence ${\bf A}_L^\dagger\subseteq{\bf A}_\Psi^\dagger$ while the inclusion ${\bf A}_\Psi^\dagger\subseteq \tilde{\bf A}^\dagger\cap {\bf A}_L={\bf A}_L^\dagger$ follows from Prop.\  \ref{prop:Lpsi}.
\end{proof}

 In Definition 5.6.1 in  (loc.\ cit.) they define the property {\it decompleting} for a tower $\Psi,$ which we are not going to recall here as it is rather technical. The cyclotomic tower over $\Qp$ is of this kind for instance. If our $\Psi$ would be decompleting, the machinery of (loc.\ cit.), in particular Theorems  5.7.3/4, adapted to the Lubin-Tate setting would imply that  all the categories at the end of section \ref{subsec:picture} are equivalent, which contradicts Remark \ref{rem:notoverconvergent}.

\noindent
Peter Schneider,\\
Universit\"{a}t M\"{u}nster,  Mathematisches Institut,\\  Einsteinstr. 62,
48291 M\"{u}nster,  Germany,\\
 http://www.uni-muenster.de/math/u/schneider/  \\
pschnei@uni-muenster.de \\ \\

\noindent
Otmar Venjakob\\
Universit\"{a}t Heidelberg,  Mathematisches Institut,\\  Im Neuenheimer Feld 288,  69120
Heidelberg,  Germany,\\
 http://www.mathi.uni-heidelberg.de/$\,\tilde{}\,$venjakob/\\
venjakob@mathi.uni-heidelberg.de

%\printendnotes[custom]


\newpage

\addcontentsline{toc}{part}{\large References}

\begin{thebibliography}{B-GAL}
%
%%\bibitem[Bas]{Bas}
%%Bass H.: \emph{Algebraic $K$-Theory}. Benjamin 1968
%%
%%\bibitem[Col]{Col}
%%Coleman R.: \emph{Division Values in Local Fields}. Invent.\ math.\ 53, 91-116 (1979)
%%
%%\bibitem[CE]{CE}
%%Chiarellotto, B., Esposito, F. \emph{A note on Fontaine theory using different Lubin-Tate groups.} Kodai Math. J. 37 (2014), no. 1, 196--211.
%%
%%\bibitem[Fo]{Fo}
%%Fontaine, J.M.: \emph{ Repr\'{e}sentations $p$-adiques des corps locaux} The Grothendieck Festschrift \textbf{2}, Birkh\"{a}user, Boston, 1991, 249--309.
%%\bibitem[FX]{FX}
%%Fourquaux L., Xie B.: \emph{Triangulable $O_F$-analytic $(\varphi_q,\Gamma)$-modules of rank $2$}. arXiv:1206.2102v3 (2013)
%%
%%
%%\bibitem[Jen]{Jen}
%%Jensen C.U.: \emph{Les Foncteurs D\'eriv\'es de $\varprojlim$ et leurs Applications en Th\'eorie des Modules}.
%%Springer Lect.\ Notes Math., vol.\ 254 (1972)
%%
%%\bibitem[Haz]{haz}
%%Hazewinkel M.: \emph{Formal groups and applications}.
%%Academic Press 1978
%%
%%\bibitem[Her1]{her98}
%%Herr, L.: \emph{Sur la cohomologie galoisienne des corps $p$-adiques. (French) [On the Galois cohomology of $p$-adic fields]} Bull. Soc. Math. France 126 (1998), no. 4, 563--600.
%%
%%
%%\bibitem[Her2]{her}
%%Herr, L.: \emph{Une approche nouvelle de la dualit\'{e} locale de Tate. (French) [A new approach to Tate's local duality]} Math. Ann. 320 (2001), no. 2, 307--337.
%\bibitem[AL]{AL}
%Abe, T.; Lazda, C.:\emph{Proper pushforwards on finite dimensional adic spaces.} arXiv preprint arXiv:2009.05433 (2020).




\bibitem[Ax]{Ax}
Ax, J.:\emph{ Zeros of polynomials over local fields—The Galois action.} J. Algebra 15 (1970), 417--428.
%
%\bibitem[BS]{BS}
%B\v{a}nic\v{a} C., St\v{a}n\v{a}\c{s}il\v{a} O.: \emph{Algebraic Methods in the Global Theory of Complex Spaces}. J.\ Wiley \& Sons 1976
%
\bibitem[BV]{BV}
Bellovin, R., Venjakob, O.:\emph{Wach modules, regulator maps, and $\epsilon$-isomorphisms in families.} Int. Math. Res. Not. IMRN 2019, no. 16, 5127–5204.


\bibitem[Ben]{benois} Benois, D.: \emph{ On Iwasawa theory of crystalline representations.} Duke Math. J. 104, no. 2, 211--267  (2000)
%
\bibitem[B]{B} Berger, L.: \emph{Bloch and Kato's exponential map: three explicit formulas.} Kazuya Kato's fiftieth birthday. Doc. Math., Extra Vol., 99-129 (2003)
%
%\bibitem[Be]{Be}
%Berger L.: \emph{Limites de repr\'esentations cristallines}.  Compositio Math.\ 140, 1473-1498 (2004)

\bibitem[Be16]{Be16}
Berger, L.: \emph{ Multivariable $(\varphi,\Gamma)$-modules and locally analytic vectors.} Duke Math. J. 165 , no. 18, 3567--3595 (2016)

\bibitem[BeCo]{BeCo}
\textit{L. Berger} and \textit{P. Colmez}:\emph{Familles de repr\'{e}sentations de de Rham et monodromie $p$-adique}, Ast\'{e}risque 319, 303--337 (2008)

\bibitem[BF]{BF} Berger L., Fourquaux L.: \emph{Iwasawa theory and $F$-analytic $(\varphi,\Gamma)$-modules}. Preprint 2015
%
\bibitem[BSX]{BSX} Berger L., Schneider P., Xie B.: \emph{Rigid character groups, Lubin-Tate theory, and $(\varphi,\Gamma)$-modules}. To appear in Memoirs AMS
%%
%
%\bibitem[Ber]{Ber}
%Berkovich V.: \emph{Spectral theory and analytic geometry over non-archimedean fields}. AMS 1990
%
%\bibitem[Bey]{Bey}
%Beyer P.: \emph{Zur Serre-Dualit\"at f\"ur koh\"arente Garben auf rigid-analytischen R\"aumen}. Schriftenreihe Math.\ Institut Uni M\"unster, 3.\ Serie, Heft 20, 1-64 (1997)
%
%\bibitem[BK]{BK}
%Bloch S.,  Kato K.: \emph{$L$-functions and Tamagawa numbers of motives}. The Grothendieck Festschrift, Vol.\ I, 333-400, Progress Math., 86, Birkh\"auser Boston 1990
%
%\bibitem[BW]{BW}
%Borel, A.; Wallach, N. R.: \emph{Continuous cohomology, discrete subgroups, and representations of reductive groups.} Annals of Mathematics Studies, 94. Princeton University Press, Princeton, N.J.; University of Tokyo Press, Tokyo, 1980.
%
%\bibitem[BGR]{BGR}
%Bosch S., G\"untzer U., Remmert R.: \emph{Non-Archimedean Analyis}. Springer Grundlehren, vol.\ 261, Heidelberg 1984
%
%\bibitem[B-TG]{B-TG}
%Bourbaki N.: \emph{Topologie G\'en\'erale}. Chap.\ 1--4, 5--10. Springer 2007
%
%\bibitem[B-TVS]{B-TVS}
%Bourbaki N.: \emph{Topological Vector Spaces}. Springer 1987
%
%\bibitem[B-CA]{B-CA}
%Bourbaki N.: \emph{Commutative Algebra}. Hermann

\bibitem[BC]{BC}
Brinon,  O.; Conrad, B.: \emph{ Notes on p-adic Hodge theory.} Notes from the CMI Summer School, preprint, 2009

\bibitem[Coh]{Coh}
Cohn, P.M.: \emph{Algebra, Volume 3.}  Second edition. John Wiley \& Sons, Ltd., Chichester, 1991.
\bibitem[Ede]{Ede}
Edenfeld, V.: \emph{A pre-perfectoid approach to Robba rings.} Dissertation, M\"{u}nster 2022.


\bibitem[ChCo1]{ChCo}
Cherbonnier, F.; Colmez, P.: \emph{ Repr\'{e}sentations $p$-adiques surconvergentes.} Invent. Math. 133 , no. 3, 581--611 (1998)


%\bibitem[ChCo2]{ChCo2}
%Cherbonnier, F.; Colmez, P.: \emph{ Th\'{e}orie d'Iwasawa des repr\'{e}sentations $p$-adiques d'un corps local.}   J. Amer. Math. Soc. 12 (1999), no. 1, 241--268.
%
%\bibitem[Chi]{Chi}
%Chiarellotto B.: \emph{Duality in Rigid Analysis}. In p-adic Analysis (eds.\ Baldassarri, Bosch, Dwork), Springer Lect.\ Notes Math., vol.\ 1454, pp.\ 142-172 (1990)
%
%\bibitem[CR]{CR}
%Christol G., Robba P.: \emph{\'Equations diff\'erentielles $p$-adiques}. Hermann 1994
%
%\bibitem[Co1]{Co}
%Colmez P.: \emph{Repr\'esentations cristallines et repr\'esentations de hauteur finie}. J.\ reine angew.\ Math.\ 514, 119-143 (1999)
%

\bibitem[Co1]{Co}
Colmez P.: \emph{Repr\'{e}sentations   de $GL_2(\Qp)$ et $(\varphi,\Gamma)$-modules}. Ast\'{e}risque (2010), p. 281-509.

\bibitem[Co2]{Co2}
Colmez P.: \emph{Repr\'{e}sentations localement analytiques de $GL_2(\Qp)$ et $(\varphi,\Gamma)$-modules}. Representation Theory 20,  187–248 (2016)
%
%\bibitem[Co3]{Co3}
%Colmez P.: \emph{ $(\varphi,\Gamma)$-modules et repr\'{e}sentations du mirabolique de ${\rm GL}_2(\mathbf Q_p)$.}   Ast\'{e}risque No. 330, 61--153 (2010)
%
%\bibitem[Co4]{Co4}
%Colmez P.: \emph{Espaces de Banach de dimension finie}. J.\ Inst.\ Math.\ Jussieu 1, 331-439 (2002)
%
%\bibitem[Co5]{CoIw}
%Colmez, P.: \emph{ Th\'eorie d'Iwasawa des repr\'esentations de de Rham d'un corps local}. Ann.\ Math.\ 148 (2), 485-571 (1998)
%
\bibitem[CoNi]{CoNi}
Colmez, P. and Niziol, W.:\emph{ Syntomic complexes and $p$-adic nearby cycles.} Invent. Math. 208, no. 1, 1--108  (2017)
%
%\bibitem[EGA]{EGA}
%Grothendieck A., Dieudonn\'e J.: \emph{\'El\'ement de G\'eom\'etrie Alg\'ebrique I}. Springer 1971
%
%\bibitem[Ei]{Ei}
%Eisenbud, D.: \emph{Commutative algebra. With a view toward algebraic geometry.} Graduate Texts in Mathematics, 150. Springer-Verlag, New York, 1995.
%
%
%\bibitem[Eme]{Eme}
%Emerton M.: \emph{Locally analytic vectors in representations of locally $p$-adic analytic groups}. Memoirs AMS 248, 1175 (2017)
%
\bibitem[FF]{FF} Fargues, L. and Fontaine, J.-M.: \emph{ Courbes et Fibr\'{e}s Vectoriels en Th\'{e}orie de Hodge $p$-adique.}   Ast\'{e}risque, 406, Soc. Math. France, Paris, 2018.
%
%\bibitem[F1]{F1} Fontaine, J.-M.: \emph{ Modules galoisiens, modules filtr\'{e}s et anneaux de Barsotti-Tate.} Journ\'{e}es de G\'{e}om\'{e}trie Alg\'{e}brique de Rennes. (Rennes, 1978), Vol. III, pp. 3--80, Ast\'{e}risque, 65, Soc. Math. France, Paris, 1979.
%
%\bibitem[F2]{F}   Fontaine J.-M.: \emph{ Sur certains types de repr\'{e}sentations $p$-adiques du groupe de Galois d'un corps local; construction d'un anneau de Barsotti-Tate.}    Ann. of Math. (2) 115, no. 3, 529--577 (1982)

\bibitem[Fo]{Fo}
Fontaine J.-M.: \emph{R\'epresentations $p$-adiques des corps locaux}. In ``The Grothendieck Festschrift'', vol.\ II, 249-309, Birkh\"auser 1990

\bibitem[FO]{FO}   Fontaine J.-M., Ouyang Y.: \emph{ Theory of $p$-adic Galois representations.}     preprint.

\bibitem[FX]{FX}
Fourquaux L., Xie B.: \emph{Triangulable $O_F$-analytic $(\varphi_q,\Gamma)$-modules of rank $2$}.  Algebra \& Number Theory 7 (10), 2545-2592 (2013)

%\bibitem[FK]{FK}
%Fukaya T., Kato K.: \emph{A formulation of conjectures on $p$-adic zeta functions in non-commutative Iwasawa theory}. Proc.\ St.\  Petersburg Math.\ Soc., Vol.\ XII, AMS Transl.\ Ser.\ 2, vol.\ 219, 1-86 (2006)
%
%\bibitem[G]{G}
%Gruson, L.: \emph{ Th\'{e}orie de Fredholm $p$-adique.}  Bull. Soc. Math. France 94 (1966), 67--95.
%
%\bibitem[Haz]{Haz}
%Hazewinkel M.: \emph{Formal Groups and Applications}. Academic Press 1978

\bibitem[Her98]{Her1}
Herr, L.: \emph{Sur la cohomologie galoisienne des corps $p$-adiques.}   Bull. Soc. Math. France 126 (1998), no. 4, 563--600.

\bibitem[Ked05]{Ked05}
Kedlaya, K.: \emph{Slope filtrations revisited.} Doc. Math. 10 (2005), 447--525.

\bibitem[Ked08]{Ked08}
Kedlaya, K.: \emph{Slope filtrations for relative Frobenius.} Repr\'{e}sentations $p$-adiques de groupes $p$-adiques. I. Repr\'{e}sentations galoisiennes et $(\varphi,\Gamma)$-modules. Ast\'{e}risque No. 319 (2008), 259--301.

%\bibitem[Ked]{Ked}
%Kedlaya K.: \emph{Some slope theory for multivariate Robba rings}. arXiv:1311.7468v1

\bibitem[Ked15]{KedNew}
Kedlaya, K.: \emph{New methods for $(\varphi,\Gamma)$-modules.} Res. Math. Sci. 2 (2015), Art. 20, 31 pp.

\bibitem[KP]{KP}
Kedlaya K., Pottharst J.: \emph{On categories of $(\varphi,\Gamma)$-modules.} Algebraic geometry: Salt Lake City 2015, 281--304, Proc. Sympos. Pure Math., 97.2, Amer. Math. Soc., Providence, RI, 2018.

\bibitem[KPX]{KPX}
Kedlaya K., Pottharst J. and Xiao L.: \emph{Cohomology of arithmetic families of $(\varphi,\Gamma)$-modules}. (Zitate aus arXiv:1203.5718v1! Aktualisieren!?)  J. Amer. Math. Soc. 27 (2014), no. 4, 1043--1115.

\bibitem[KLI]{KLI}
Kedlaya K., Liu, R.: \emph{Relative $p$-adic Hodge theorie: foundations}.  Ast\'{e}risque No. 371 (2015), 239 pp.

\bibitem[KLII]{KLII}
Kedlaya K., Liu, R.: \emph{Relative $p$-adic Hodge theorie II}.

%\bibitem[KLIII]{KLIII}
%Kedlaya K., Liu, R.: \emph{Finiteness of cohomology of local systems on rigid analytic spaces,} arXiv:1611.06930v1, 2016.
%
%
%\bibitem[Kis]{Kis}
%Kisin M.: \emph{Crystalline representations and $F$-crystals}. Algebraic geometry and number theory, Progress Math., vol.\ 253, pp.\ 459-496, Birkh\"auser 2006


\bibitem[KR]{KR} Kisin M., Ren W.: \emph{ Galois representations and Lubin-Tate groups}. Doc. Math., vol.\ 14 , 441-461 (2009)


\bibitem[Kl]{Kl}
Kley, M.:\emph{Perfekte $(\varphi, \Gamma)$-Moduln.} Masterarbeit, M\"{u}nster 2016


\bibitem[KV]{KV}
Kupferer, B., Venjakob,O.: \emph{{Herr-complexes in the Lubin-Tate setting}},
   {2020}.
%
\bibitem[Ku]{Ku}
Kupferer, B.: \emph{{Two ways to compute Galois Cohomology using Lubin-Tate
  {$(\varphi,\Gamma)$}-Modules, a Reciprpcity Law and a Regulator Map}},
  {Ruprecht-Karls-Universit\"{a}t Heidelberg -
  \url{https://www.mathi.uni-heidelberg.de/~otmar/doktorarbeiten/Dissertation
  Benjamin Kupferer.pdf}}, {2020}.
%

%\bibitem[Lan]{Lan}
%Lang S.: \emph{Cyclotomic Fields}. Springer 1978
%
%\bibitem[Lau]{Lau}
%Laudal O.: \emph{ Projective systems and valuation theory}.  Matematisk Seminar
%Universitetet i Oslo, Nr. 5, April 1965
%
%\bibitem[Laz1]{Laz}
%Lazard M.: \emph{Les z\'eros des fonctions analytiques d'une variable sur un corps valu\'e complet}. Publ.\ Math.\ IHES 14, 47-75 (1962)
%
\bibitem[Laz]{La}
Lazard, M.: \emph{ Groupes analytiques $p$-adiques.}   Inst. Hautes \'{E}tudes Sci. Publ. Math. No. 26  389--603 (1965)
%
%\bibitem[Li]{Li}
%Liu,  R.: \emph{ Cohomology and duality for $(\varphi,\Gamma)$-modules over the Robba ring.} Int. Math. Res. Not. IMRN  no. 3, Art. ID rnm150, 32 pp. (2008)
%
\bibitem[LLZ11]{LLZ11}
Lei, A., Loeffler, D.; Zerbes, S. L.:\emph{
Coleman maps and the p-adic regulator. } 
Algebra Number Theory 5 (2011), no. 8, 1095–1131. 
%
\bibitem[LVZ15]{LVZ15}
Loeffler D., Venjakob O., Zerbes S.\ L.: \emph{Local epsilon isomorphisms}. Kyoto J.\ Math.\ \textbf{55} (2015), no. 1, 63--127.
%
%\bibitem[LZ]{LZ}
%Loeffler, D., Zerbes, S.\ L.: \emph{ Iwasawa theory and $p$-adic $L$-functions over $\mathbb{Z}_p^2$-extensions.} Int. J. Number Theory 10 (2014), no. 8, 2045--2095.
%
%\bibitem{Lu}{Lu0}
%{Lubkin}, S.: \emph{{Cohomology of completions.}}, vol.~42, Elsevier, Amsterdam,
%  1980.
%
\bibitem[MSVW]{MSVW} 
Milan Malcic,M., Steingart,R., Venjakob,O. and Witzelsperger,M.:\emph{
$\epsilon$-Isomorphisms for rank one Lubin-Tate $(\varphi,\Gamma)$-modules over the Robba ring,} preprint (2023).



%
%\bibitem[Mal]{Mal}
%Malcic M.: \emph{A relative trace map and its compatibility with Serre duality in rigid analytic
%geometry }. PhD thesis, Heidelberg, 2022, in preparation
%
%\bibitem[Mat]{Mat}
%Matsumura H.: \emph{Commutative ring theory}. Cambridge Univ.\ Press 1986
%
%\bibitem[MCR]{MCR}
%McConnell J.C., Robson J.C.: \emph{Noncommutative Noetherian Rings}. AMS 2001
%
%\bibitem[MS]{MS}
%Morita Y., Schickhof W.: \emph{Duality of projective limit spaces and inductive limit spaces over a nonspherically complete nonarchimedean field}. T\^{o}hoku Math.\ J.\ 38(3), 387-397 (1986)
%

\bibitem[Na14a]{Na}
 Nakamura, Kentaro. Iwasawa theory of de Rham $(\varphi,\Gamma)$-modules over the Robba ring. J. Inst. Math. Jussieu 13 (2014), no. 1, 65--118.


\bibitem[Na17a]{NaANT}
Nakamura, Kentaro. A generalization of Kato's local $\varepsilon$-conjecture for $(\varphi,\Gamma)$-modules over the Robba ring. Algebra Number Theory 11 (2017), no. 2, 319--404.

\bibitem[Na17b]{NaRk}
Nakamura, Kentaro. Local $\varepsilon$-isomorphisms for rank two $p$-adic representations of ${\rm Gal}(\overline{\Bbb Q}_p/\Bbb Q_p)$ and a functional equation of Kato's Euler system. Camb. J. Math. 5 (2017), no. 3, 281--368.

%\bibitem[Ne]{Ne}
%Nekov\'{a}{$\mathrm{\check{r}}$}, J.: \emph{ Selmer complexes.} Ast\'{e}risque No. 310 (2006)
%
%\bibitem[Ni]{niziol}
%Niziol, W.: \emph{ Cohomology of crystalline representations.} Duke Math. J. \textbf{71}, no. 3, 747--791  (1993)

\bibitem[NSW]{NSW}
Neukirch J., Schmidt A., Wingberg K.: \emph{Cohomology of Number Fields}. $2^{nd}$ Ed., Springer 2008

%\bibitem[PZ]{PZ}
%Pal A., Zabradi G.: \emph{Cohomology and overconvergence for representations of powers of Galois groups}. arXiv:1705.03786v4 2019
%
%\bibitem[PGS]{PGS}
%Perez-Garcia C., Schikhof W.H.: \emph{ Locally Convex Spaces over Non-Archimedean Valued Fields}. Cambridge Univ.\ Press 2010
%%
%%\bibitem[P]{P}
%%Perrin-Riou B.: \emph{Theorie d'Iwasawa $p$-adique locale et globale}. Invent.\ Math.\ 99, 247-292 (1990)
%
\bibitem[Poy]{Poy}
Poyeton, L.:\emph{  A note on $F$-analytic $B$-pairs},
    eprint arXiv:2011.04900, (2020)
%
%\bibitem[Pro]{Pro}
%Prosmans F.: \emph{Derived Projective Limits of Topological Abelian Groups}. J.\ Functional Analysis 162, 135-177 (1999)
%
%\bibitem[PR]{PR}
%Perrin-Riou B.: \emph{ Th\'{e}orie d'Iwasawa des repr\'{e}sentations $p$-adiques sur un corps local.}  Invent.\ Math.\ 115, no.\ 1, 81--161   (1994)
%
%\bibitem[NFA]{NFA}
%Schneider P.: \emph{Nonarchimedean Functional Analysis}. Springer 2002

%\bibitem[Po]{Po}
%Pottharst, J.:\emph{ Analytic families of finite-slope Selmer groups.} Algebra Number Theory 7, no. 7, 1571--1612 (2013)
%
%\bibitem[Sa]{Sa}Saavedra Rivano, N.:\emph{ Cat\'{e}gories Tannakiennes.} Lecture Notes in Mathematics, Vol. 265. Springer-Verlag, Berlin-New York, 1972.
%
%\bibitem[Sc1]{Sc1}
%Schmidt, T.: \emph{Auslander regularity of $p$-adic distribution algebras.} Represent. Theory 12 (2008), 37--57.
%
%\bibitem[Sc2]{Sc}
%Schmidt T.: \emph{On locally analytic Beilinson-Bernstein localization and the canonical dimension},  Math. Z. 275 (2013), no. 3-4, 793--833
%
%\bibitem[Sch]{Sch}
%Schneider P.: \emph{Points of rigid analytic varieties}. J.\ reine angew.\ Math.\ 434, 127-157 (1993)
%
%\bibitem[NFA]{NFA}
%Schneider P.: \emph{Nonarchimedean Functional Analysis}. Springer 2002

%\bibitem[TdA]{TdA}
%Schneider P.: \emph{Die Theorie des Anstieges}. Course, M\"unster 2006/07

\bibitem[pLG]{pLG}
Schneider P.: \emph{$p$-Adic Lie Groups}. Springer Grundlehren math.\ Wissenschaften, vol.\ 344. Springer 2011

%\bibitem[App]{App}
%Schneider P.: \emph{Robba rings for compact $p$-adic Lie groups}.
%Appendix to the paper ``Generalized Robba rings'' by G.\ Zabradi, Israel J.\ Math.\ 191, 856 - 887 (2012)

\bibitem[GAL]{GAL}
Schneider P.: \emph{Galois representations and $(\varphi,\Gamma)$-modules}. Cambridge studies in advanced mathematics, vol.\ 164. Cambridge Univ.\ Press 2017

%\bibitem[ST1]{ST1}
%Schneider P., Teitelbaum J.: \emph{Locally analytic distributions and $p$-adic representation theory, with applications to $GL_2$}. J.\ AMS 15, 443-468 (2001)
%
%\bibitem[ST2]{ST2}
%Schneider P., Teitelbaum J.: \emph{$p$-adic Fourier theory}. Documenta Math.\ 6, 447-481 (2001)
%
%\bibitem[ST]{ST}
%Schneider P., Teitelbaum J.: \emph{Algebras of $p$-adic distributions and admissible representations}. Invent.\ math.\ 153, 145-196 (2003)
%
%\bibitem[ST3]{ST3}
%Schneider P., Teitelbaum J.: \emph{Duality for admissible locally analytic representations}. Representation Theory 9, 297-326 (2005)
%
%\bibitem[ST4]{ST4}
%Schneider P., Teitelbaum J.: \emph{Banach-Hecke algebras and $p$-adic Galois representations}. Documenta Math., The Book Series 4 (J.\ Coates' Sixtieth Birthday), pp.\ 631 - 684 (2006)
%
%\bibitem[SV]{SV}
%Schneider P., Venjakob O.: \emph{A splitting for $K_1$ of completed group rings}. Comment.\ Math.\ Helv.\ 88, 613-642 (2013)


\bibitem[SV15]{SV15}
Schneider P., Venjakob O.: \emph{Coates-Wiles homomorphisms and Iwasawa cohomology for Lubin-Tate extensions}.  (2015)


\bibitem[SV23]{SV23}
Schneider P., Venjakob O.: \emph{{Reciprocity laws for $(\varphi_L,\Gamma_L)$-modules over Lubin-Tate extensions}}.  (2023)



%\bibitem[Se0]{Se0}
%Serre J.-P.: \emph{Abelian $l$-Adic Representations and Elliptic
%Curves}. W.A.\ Benjamin 1968
%
%\bibitem[SP]{stacks-project}
%The {Stacks project authors}, \emph{{The Stacks project}},
%  {\url{https://stacks.math.columbia.edu}}, {2018}.
%
\bibitem[Ste]{Ste}
Steingart, R.: \emph{Frobeniusregularisierung und
Limites $L$-kristalliner
Darstellungen.} Master thesis, Heidelberg 2019

\bibitem[Ste1]{Stein}
Thomas, O.: \emph{Analytic cohomology of families of L-analytic
Lubin-Tate $(\varphi_L, \Gamma_L)$-modules.} PhD thesis, 2022, Heidelberg

\bibitem[Ta]{Ta}
Tate, J. T.: \emph{ $p$-divisible groups.} 1967 Proc. Conf. Local Fields (Driebergen, 1966) pp. 158--183 Springer, Berlin


\bibitem[V13]{V-Kato}
Otmar Venjakob, \emph{On {K}ato's local {$\epsilon$}-isomorphism conjecture for
  rank-one {I}wasawa modules}, Algebra Number Theory \textbf{7} (2013), no.~10,
  2369--2416.
%\bibitem[Th]{Th}
%Thomas, O.: \emph{On Analytic and Iwasawa Cohomology.} PhD thesis, 2019, Heidelberg
%
%
%
%
%\bibitem[Th1]{Th1}
%Thomas, O.: \emph{Cohomology of topologised monoids}, preprint 2019
%
%%\bibitem[Sch]{sch}
%%Scholl, A.~J.: \emph{Higher fields of norms and $(\phi,\Gamma)$-modules.} Doc. Math. 2006, Extra Vol., 685-709
%%
%%\bibitem[V]{ven}
%%Venjakob, O. \emph{On the Iwasawa theory of $p$-adic Lie extensions.} Compositio Math. 138 (2003), no. 1, 1-54.
%
%\bibitem[vdP]{vdP}
%van der Put M.: \emph{Serre duality for rigid analytic spaces}. Indagationes Math.\ , N.S.\ 3(2), 219-235 (1992)
%
%\bibitem[vR]{R}
%van Rooij, A. C. M.:\emph{ Non-Archimedean functional analysis.} Monographs and Textbooks in Pure and Applied Math., 51. Marcel Dekker, Inc., New York, 1978.
%
\bibitem[W]{W}
Witzelsperger, M.: \emph{Eine Kategorien\"{a}quivalenz zwischen Darstellungen und $(\varphi,\Gamma)$-Moduln \"{u}ber dem Robba-Ring.} Master thesis, Heidelberg 2020

\end{thebibliography}

\newpage
\begin{thebibliography}{B-GAL}


\bibitem[Scho]{Scho}
Scholze, P.: \emph{ $p$-adic Hodge theory for rigid-analytic varieties.} Forum Math. Pi 1 (2013)


\end{thebibliography}
\end{document}